\documentclass[leqno]{aadmbook}
\usepackage{amsthm,amsfonts,amsmath,amssymb}
\usepackage{graphicx,cite,color,hyperref,booktabs}
\usepackage[T1]{fontenc}

\newtheorem{theorem}{Theorem}
\newtheorem{lemma}{Lemma}
\newtheorem*{conjecture}{Conjecture AAA}

\newcommand*\diff{\mathop{}\!\mathrm{d}}

\usepackage{float}
\usepackage{tikz}
\usetikzlibrary{arrows}
\usetikzlibrary{arrows.meta}
\usetikzlibrary{chains}
\usetikzlibrary{positioning}
\usetikzlibrary{automata,positioning,calc}
\usetikzlibrary{decorations}
\usetikzlibrary{decorations.shapes}
\usetikzlibrary{decorations.markings}
\tikzset{
    edge/.style={-{Latex[scale=1.7]}},
    dedge/.style={{Latex[scale=1.7]}-{Latex[scale=1.7]}},
}

\begin{document}
\thispagestyle{plain}

\begin{center}
{\large \sc  Applicable Analysis and Discrete Mathematics}

{\small available online at  http:/$\!$/pefmath.etf.rs }
\end{center}

\noindent{\small{\sc  Appl. Anal. Discrete Math.\ }{\bf x} (xxxx),xxx--xxx.}\\
\noindent{\scriptsize doi:10.2298/AADMxxxxxxxx}

\vspace{5cc}
\begin{center}

{\large\bf FINDING COUNTEREXAMPLES FOR A CONJECTURE OF AKBARI, ALAZEMI AND AN\DJ{}ELI\'C
\rule{0mm}{6mm}\renewcommand{\thefootnote}{}
\footnotetext{

\rule{2.4mm}{0mm}Keywords and Phrases. Graph, Energy, Maximum vertex degree, Matching number.}}

\vspace{1cc}
{\large\it  \DJ or\dj e Stevanovi\'c, Ivan Damnjanovi\'c and Dragan Stevanovi\' c\footnote{%
\scriptsize The third author is supported by the Serbian Ministry of Education, Science and Technological Development 
through the Mathematical Institute of SASA, and by the project F-159 of the Serbian Academy of Sciences and Arts.}}

\vspace{1cc}
\parbox{24cc}{{\small
For a graph~$G$,
its energy $\mathcal{E}(G)$ is the sum of absolute values of the eigenvalues of its adjacency matrix,
the matching number $\mu(G)$ is the number of edges in a maximum matching of~$G$,
while $\Delta$ is the maximum vertex degree of~$G$.
Akbari, Alazemi and An\dj{}eli\'c in
{\it Upper Bounds on the Energy of Graphs in Terms of Matching Number.} 
Appl. Anal. Discrete Math. (2021), doi:10.2298/AADM201227016A,
proved that $\mathcal{E}(G)\leq 2\mu(G)\sqrt{\Delta}$ holds when $G$ is connected and $\Delta\geq 6$,
and conjectured that the same inequality is also valid when $2\leq\Delta\leq5$.
Here we first computationally enumerate small counterexamples for this conjecture and 
then provide two infinite families of counterexamples.
}}
\end{center}


\vspace{1.5cc}
\begin{center}
{\bf 1. INTRODUCTION}
\end{center}

Developments in both computer hardware and the methods of scientific computing over the last four decades 
have had a profound impact on studies in graph theory,
by enabling researchers to shape their intuition and test their conjectures on smaller graphs 
before investing their time in an effort to prove desired results theoretically.

Among the earliest examples of software aimed to help graph theorists in their research we can select 
the expert system GRAPH,
developed by Drago\v s Cvetkovi\'c and his collaborators at the Faculty of Electrical Engineering of the University of Belgrade
from 1980 and 1984~\cite{graph2,graph3,graph1},
and the package {\tt nauty},
developed by Brendan McKay at the Australian National University since 1981~\cite{nauty1}.
The system GRAPH represented an integrated environment in which 
a researcher could enter, modify, store and visualize graphs and compute their invariants,
without a need for additional programming.
It gathered a considerable following among researchers in spectral graph theory,
being mentioned in 92 papers published from 1982--2001~\cite{graph4}.

On the other hand, the package {\tt nauty} is still being developed (check \cite{nautyCurrent} for the latest release).
While its main purpose is computation of automorphism groups of graphs,
it contains a versatile and efficient generator {\tt geng} of graphs,
which is often used in performing exhaustive computations on small graphs.
As a matter of fact,
sets of graphs generated by {\tt geng} represent a usual starting point
for the use of a more recent Java framework {\tt graph6java}~\cite{graph6java}
for programmatic testing of conjectures in graph theory.

A substantial paradigm shift occurred at the end of 1990s when it was realized that 
many results in graph theory may be modeled as optimization problems.
The program AutoGraphiX~\cite{agx1,agx2,agx3} was written with that viewpoint in mind,
utilizing a well known variable neighborhood search metaheuristic~\cite{vns}
to find graphs that either minimize or maximize an expression built up from various graph invariants.
It achieved quite some popularity in graph theoretical community
due to the fact that its optimization results were often stated in the form of conjectures 
in a lengthy ``Variable neighborhood search for extremal graphs'' series of papers 
by Pierre Hansen and his numerous coauthors,
and some of these conjectures are still being resolved today.

The recent breakthrough in machine learning algorithms may soon lead to yet another paradigm shift.
Adam Wagner~\cite{wagner} has recently used a particular reinforcement learning algorithm, the so-called deep cross-entropy method,
to construct explicit counterexamples to several conjectures in extremal combinatorics and graph theory.
Without any prior knowledge about the problem,
the learning agent here plays the game in which 
it repeatedly constructs sets of graphs that are then scored in accordance with the conjectured inequality.
The scoring function is unknown to the agent,
yet by keeping a certain percentage of top performing constructions between the iterations,
the agent slowly learns the moves which improve its performance in the game,
thus potentially leading it to the counterexamples.

Our goal here is to illustrate the use of reinforcement learning and exhaustive search
by employing them in succession to find counterexamples 
to a recent conjecture of Akbari, Alazemi and An\dj{}eli\'c \cite{aaa}.
For a graph $G$ with the vertex set $V=\{v_1,\dots,v_N\}$ and the edge set $E$,
the adjacency matrix of $G$ is the $N\times N$ matrix $A$ 
such that $A_{ij}=1$ if $v_iv_j\in E$, and $A_{ij}=0$ otherwise.
The energy of $G$ is defined as $\mathcal{E}(G)=\sum_{i=1}^N |\lambda_i|$,
where $\lambda_1\geq\cdots\geq\lambda_N$ are the eigenvalues of~$A$.
The matching number $\mu(G)$ is the number of edges in a maximum matching of $G$,
while $\Delta$ is the maximum vertex degree in $G$.
Akbari, Alazemi and An\dj{}eli\'c~\cite[Theorem~18]{aaa} proved the following theorem.
\begin{theorem}
Let $G$ be a connected graph with the maximum vertex degree $\Delta\geq 6$.
Then $\mathcal{E}(G)\leq2\mu(G)\sqrt{\Delta}$.
\end{theorem}
Although the bound $\Delta\geq 6$ was extensively used in the proof of this theorem,
the authors nevertheless posed the following conjecture as~\cite[Conjecture 23]{aaa},
which we aptly rename here.
\begin{conjecture}
$\mathcal{E}(G)\leq 2\mu(G)\sqrt{\Delta}$ holds for any (connected) graph $G\not\cong C_3, C_5, C_7$ 
with the largest vertex degree~$\Delta\in\{2,3,4,5\}$.
\end{conjecture}
In the next section we first apply Wagner's reinforcement learning approach from~\cite{wagner}
to learn that the potential counterexamples are likely to have $\Delta=3$.
In Section~3 we then apply exhaustive search to enumerate counterexamples
among connected graphs with $\Delta\leq 3$ and at most 19 vertices.
Finally, Section~4 shows that two families of graphs, inspired by the structure of the majority of small counterexamples,
represent infinite families of counterexamples for Conjecture~AAA.

\bigskip
\begin{center}
{\bf 2. SEARCHING FOR COUNTEREXAMPLES WITH REINFORCEMENT LEARNING}
\end{center}

Adam Wagner made his Python code,
which employed reinforcement learning to find counterexamples in~\cite{wagner},
publicly available at \cite{wagnerCode}.
The code requires some initial juggling to set up appropriate versions of required Python libraries.
Once that is overcome, 
the user only needs to modify the {\tt calcScore} function in the code 
in order to reuse it in a different context.
The learning agent calls this function to find about the performance of a particular graph,
and in our case, 
the score returned is $\mathcal{E}(G)-2\mu(G)\sqrt{\Delta}$ if the graph is connected,
and $-\infty$ if the graph is not connected (so as to refrain the learning agent from spending time on disconnected graphs).

\begin{figure}[h!]
\begin{center}
\includegraphics[width=0.6\textwidth]{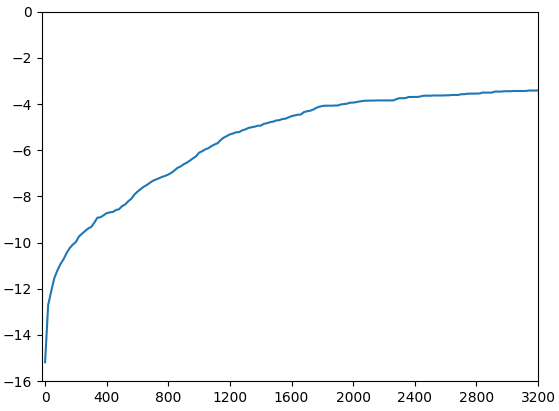}
\end{center}
\vspace{-12pt}
\caption{The mean value of $\mathcal{E}(G)-2\mu(G)\sqrt{\Delta}$ 
              for the 10\% of best performing graphs in each generation of reinforcement learning.}
\label{fig-rl-1}
\vspace{12pt}
\end{figure}

\begin{figure}[h!]
\begin{center}
\begin{tabular}{cccc}
\includegraphics[width=0.21\textwidth]{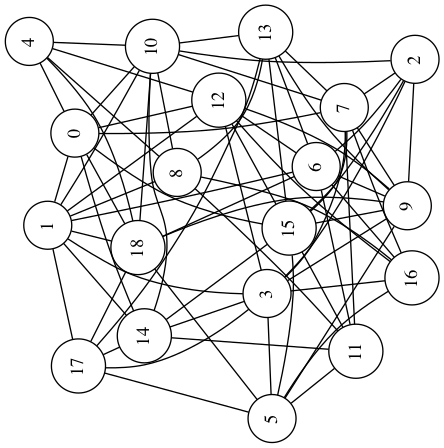} &
\includegraphics[width=0.21\textwidth]{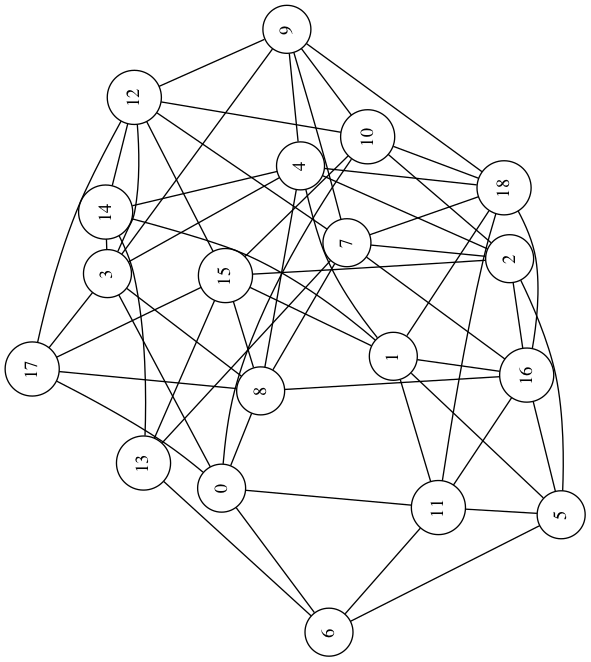} &
\includegraphics[width=0.21\textwidth]{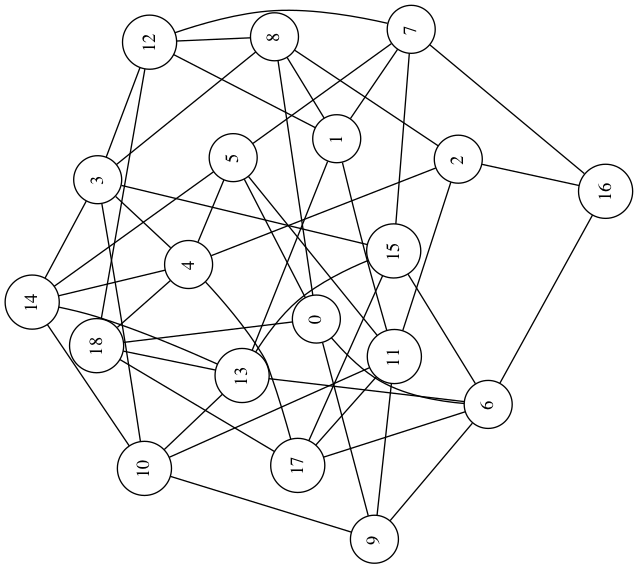} &
\includegraphics[width=0.21\textwidth]{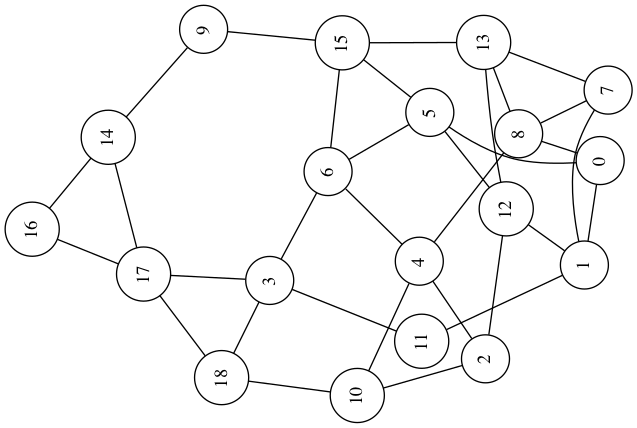} \\
200 generations & 400 generations & 600 generations & 800 generations \\[12pt]
\includegraphics[width=0.21\textwidth]{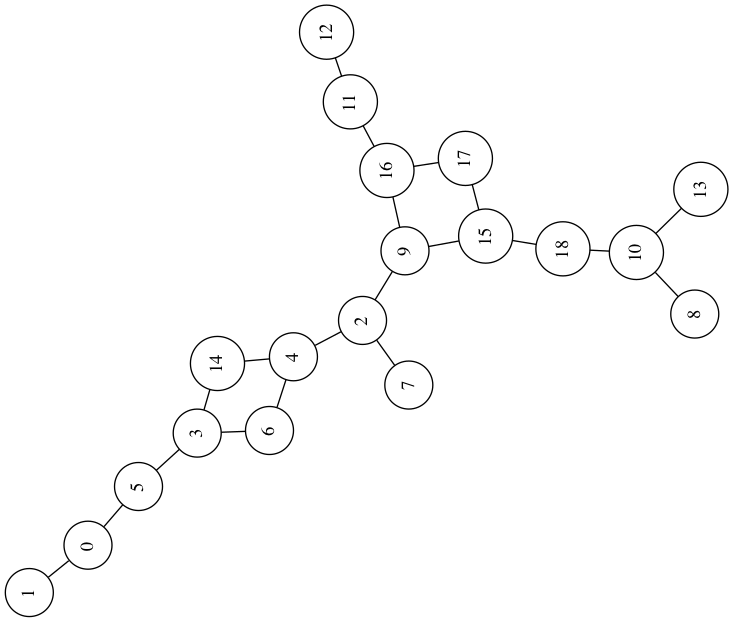} &
\includegraphics[width=0.21\textwidth]{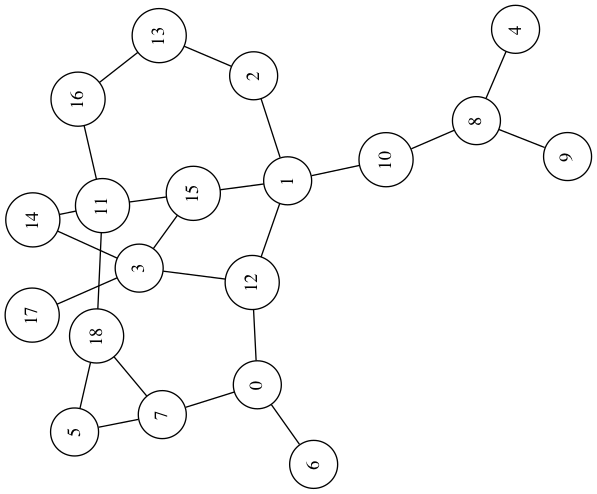} &
\includegraphics[width=0.21\textwidth]{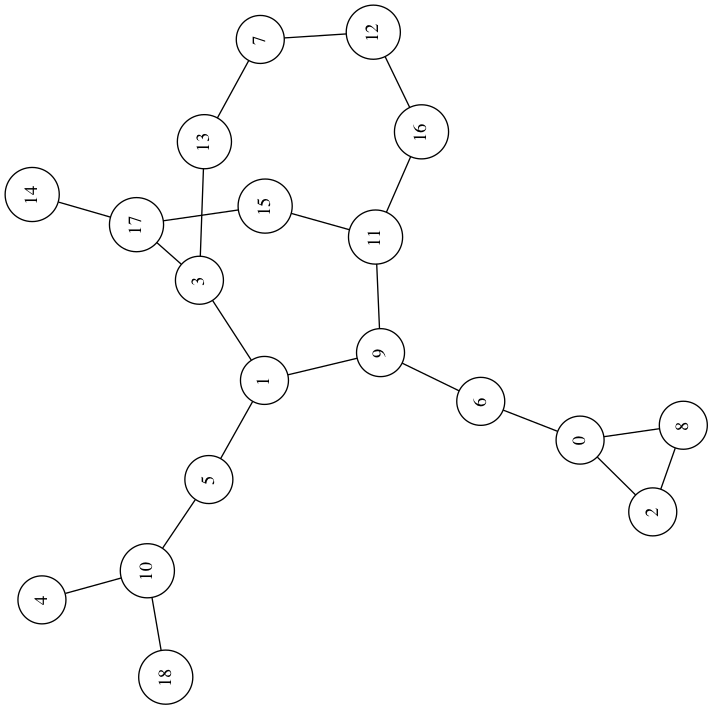} &
\includegraphics[width=0.21\textwidth]{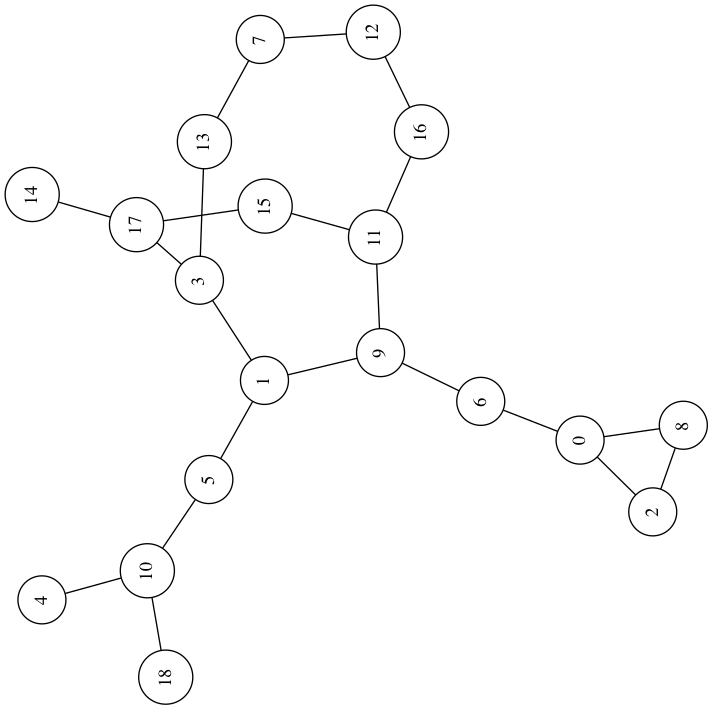} \\
1,000 generations & 1,200 generations & 1,400 generations & 1,600 generations \\[12pt]
\includegraphics[width=0.21\textwidth]{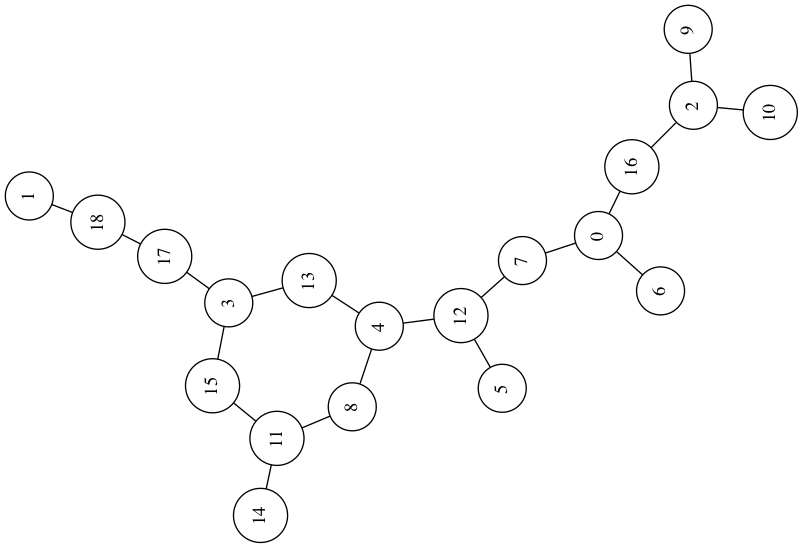} &
\includegraphics[width=0.21\textwidth]{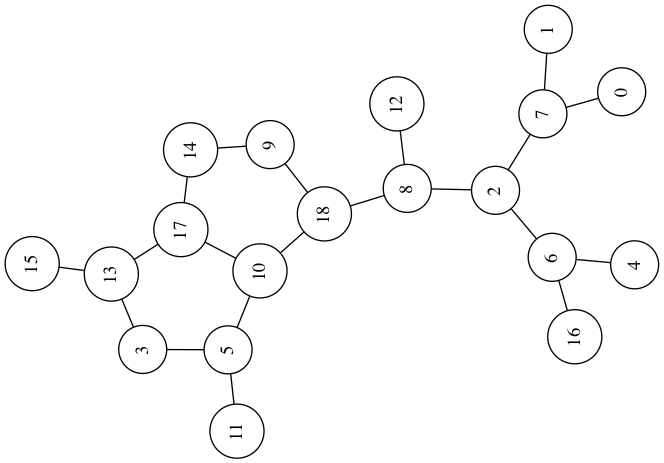} &
\multicolumn{2}{c}{\includegraphics[width=0.42\textwidth]{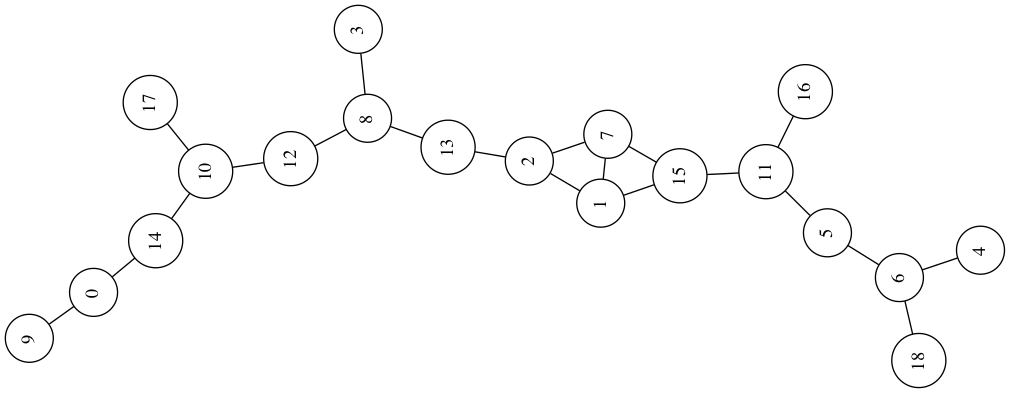}} \\
1,800 generations & 2,000 generations & 
\multicolumn{2}{c}{3,200 generations}
\end{tabular}
\end{center}
\vspace{-12pt}
\caption{The best performing graphs with the highest values of $\mathcal{E}(G)-2\mu(G)\sqrt{\Delta}$ after every 200 generations. 
             No improvement was found between the 1,400th and the 1,600th generations.
             Similarly, the same graph was the best performing one from the 2,000th until the 3,200th generation.}
\label{fig-rl-2}
\end{figure}

The learning process took place through generations:
each generation consisted of a thousand graphs on 19 vertices,
and the scoring took about ten seconds per generation on a standard i5 computer.
After 3,200 generations (which took about nine hours),
the best graph found by the learning agent had the score $\mathcal{E}(G)-2\mu(G)\sqrt{\Delta}\approx -1.665$,
meaning that no counterexample to Conjecture~AAA was found by it.
The diagram of evolution of the mean score value of the 10\% of best performing graphs in each generation,
shown in Fig.~\ref{fig-rl-1},
reveals that it was increasing at a relatively steady pace during the first 1,800 generations,
but that it almost stalled afterwards.
However,
the drawings of the best performing graphs found during the process,
shown in Fig.~\ref{fig-rl-2},
do confirm that a certain level of learning was achieved:
after 1,000 generations the best performing graphs started to be elongated structures,
and after 1,400 generations the maximum vertex degree in the best performing graphs was brought down to $\Delta=3$.

It is quite plausible that reinforcement learning would have found counterexamples to Conjecture~AAA
if it were restarted (possibly several times) with a modified scoring function
which would penalize not only disconnected, but also connected graphs with $\Delta>3$.
However, since several millions of graphs were already tested without finding a single counterexample,
we decided to proceed further with a more traditional way of searching for counterexamples,
armed with the suggestion from the reinforcement learning that
potential counterexamples should be sought among connected graphs with $\Delta\leq 3$.

\bigskip
\begin{center}
{\bf 3. FINDING COUNTEREXAMPLES WITH EXHAUSTIVE SEARCH}
\end{center}

After reinforcement learning suggested that 
potential counterexamples to Conjecture~AAA are most likely to be found among graphs with $\Delta\leq 3$,
we switched the strategy to an exhaustive search to actually enumerate small counterexamples
and learn more about their structure.
For each $6\leq N\leq 19$, 
we generated the set of such graphs with $N$ vertices with the command
\begin{quote}
{\tt geng -c -D3 $N$ > subcubic$N$.g6}
\end{quote}
(assuming the package {\tt nauty} \cite{nautyCurrent} is already installed).
On a standard i5 computer,
generation of the largest set of 317,558,689 graphs with 19 vertices lasted one hour,
while the size of the generated file was 9.84Gb.

We then modified the method {\tt run} from the class {\tt SubsetTemplate} in the Java framework {\tt graph6java}~\cite{graph6java}
so that it reports all the graphs from these sets 
for which $\mathcal{E}(G)>2\mu(G)\sqrt{\Delta}+10^{-7}$ (to account for possible numerical inaccuracies).
For the largest set of graphs with 19 vertices,
computations of energies and matching numbers lasted about 5.5 hours.
The numbers of counterexamples found in these sets 
for each number of vertices $6\leq N\leq 19$ are as follows:
\begin{center}
\begin{tabular}{rcccccccc}
\toprule
$N$                        && 6 & 7               & 8 & 9 & 10 & 11 & 12 \\
\# counterexamples && 1 & 1$^{(a)}$ & 1 & 3 &   2 &   5 &   3 \\
\midrule
$N$                        && 13 & 14 & 15 & 16 & 17 & 18 & 19 \\
\# counterexamples &&   4 &   8 &   5 & 13 & 16 & 19 & 36 \\
\bottomrule
\multicolumn{8}{l}{\footnotesize $^{(a)}$ This count excludes~$C_7$ which is forbidden by the conjecture.}
\end{tabular}
\end{center}
The first observation from this table is that, after the initial hesitation,
the number of counterexamples starts to grow steadily for $N\geq 16$,
so it seems likely that there are even more counterexamples on higher numbers of vertices.

\begin{figure}[ht]
\begin{center}
\includegraphics[width=0.95\textwidth]{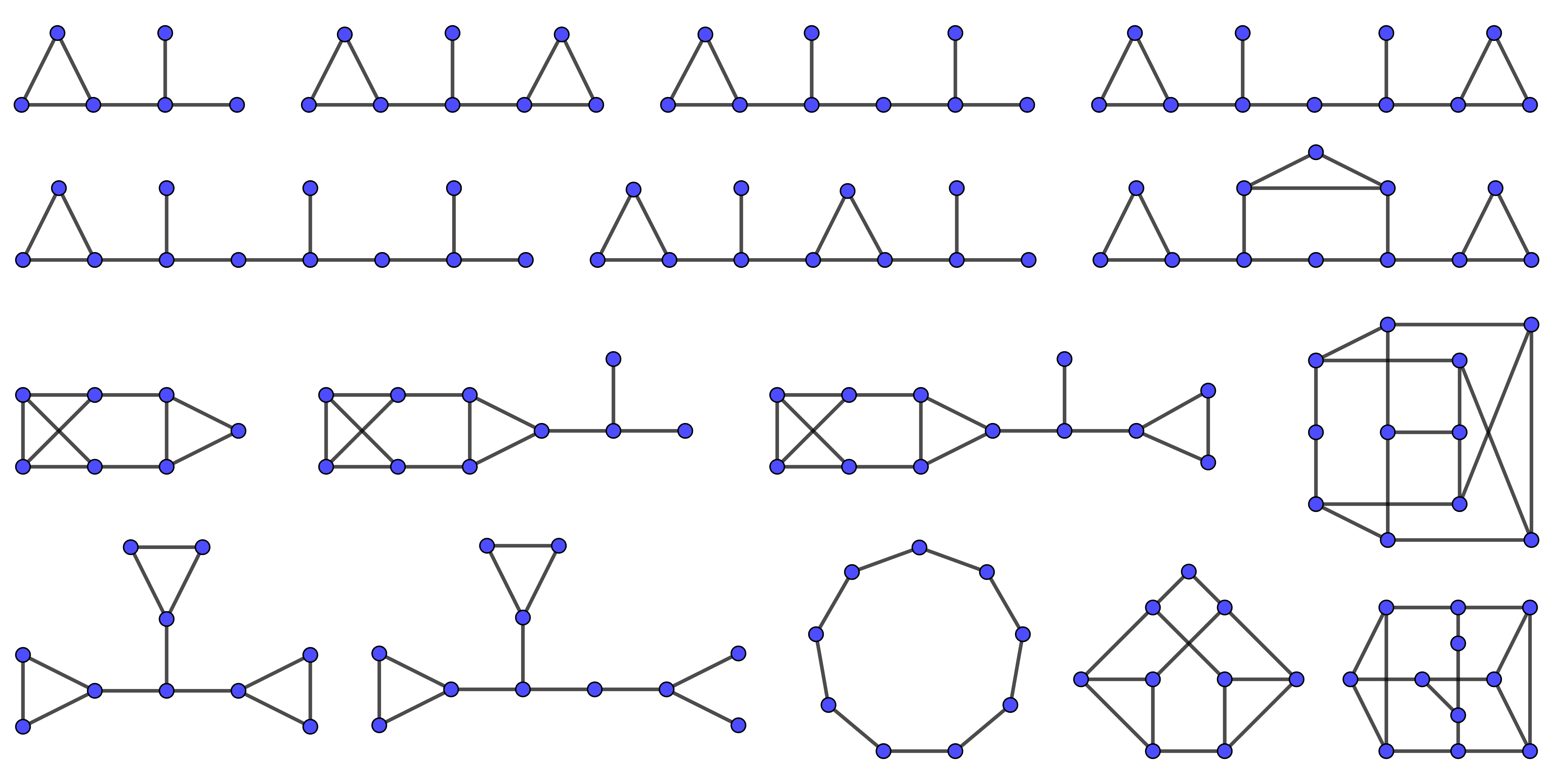}
\end{center}
\vspace{-6pt}
\caption{Counterexamples to Conjecture AAA with at most 12 vertices and $\Delta\leq 3$.}
\label{fig-ex-1}
\end{figure}

We show the drawings of counterexamples with at most 12~vertices in Fig.~\ref{fig-ex-1}.
The drawings of all counterexamples found by the exhaustive search, 
as well as the files necessary to repeat the search,
are available in~\cite{zenodo}.
We can observe from Fig.~\ref{fig-ex-1} that
the edge sets of a good number of small counterexamples may be partitioned into disjoint copies of triangles and 3-stars.
This allows for an easy ``control" of the matching number,
as any matching may contain at most one edge from any triangle or a 3-star.
Further, in most such counterexamples the copies of triangles and 3-stars are mutually connected in a tree-like manner.
This becomes an almost exclusive behaviour among larger counterexamples that have between 13 and 19 vertices.

\begin{figure}[ht]
\begin{center}
\includegraphics[scale=0.55]{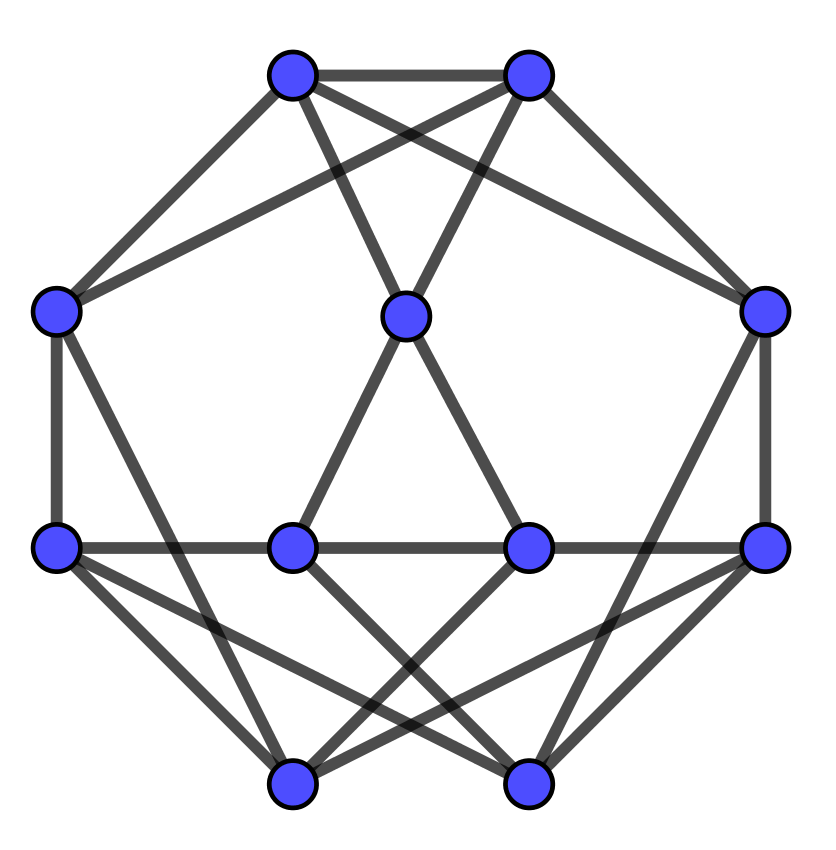}
\end{center}
\vspace{-6pt}
\caption{A counterexample to Conjecture AAA with $\Delta=4$ on 11 vertices.}
\label{fig-ex-2}
\end{figure}

For the sake of completeness,
we ran the exhaustive search on connected graphs with $\Delta\leq 4$ and $\Delta\leq 5$ as well.
We generated sets of connected graphs with $\Delta\leq 4$ and at most 14 vertices,
with the largest set containing 748 million graphs on 14 vertices.
Among these, there is a a single counterexample with $\Delta=4$ on 11 vertices,
which is shown in Fig.~\ref{fig-ex-2}.
We also generated sets of connected graphs with $\Delta\leq 5$ and at most 12 vertices,
with the largest set containing 471 million graphs on 12 vertices,
but we have not found any counterexample with $\Delta=5$ among these graphs.

\begin{center}
{\bf 4. TWO INFINITE FAMILIES OF COUNTEREXAMPLES}
\end{center}

\begin{figure}[ht]
\begin{center}
\includegraphics[width=0.75\textwidth]{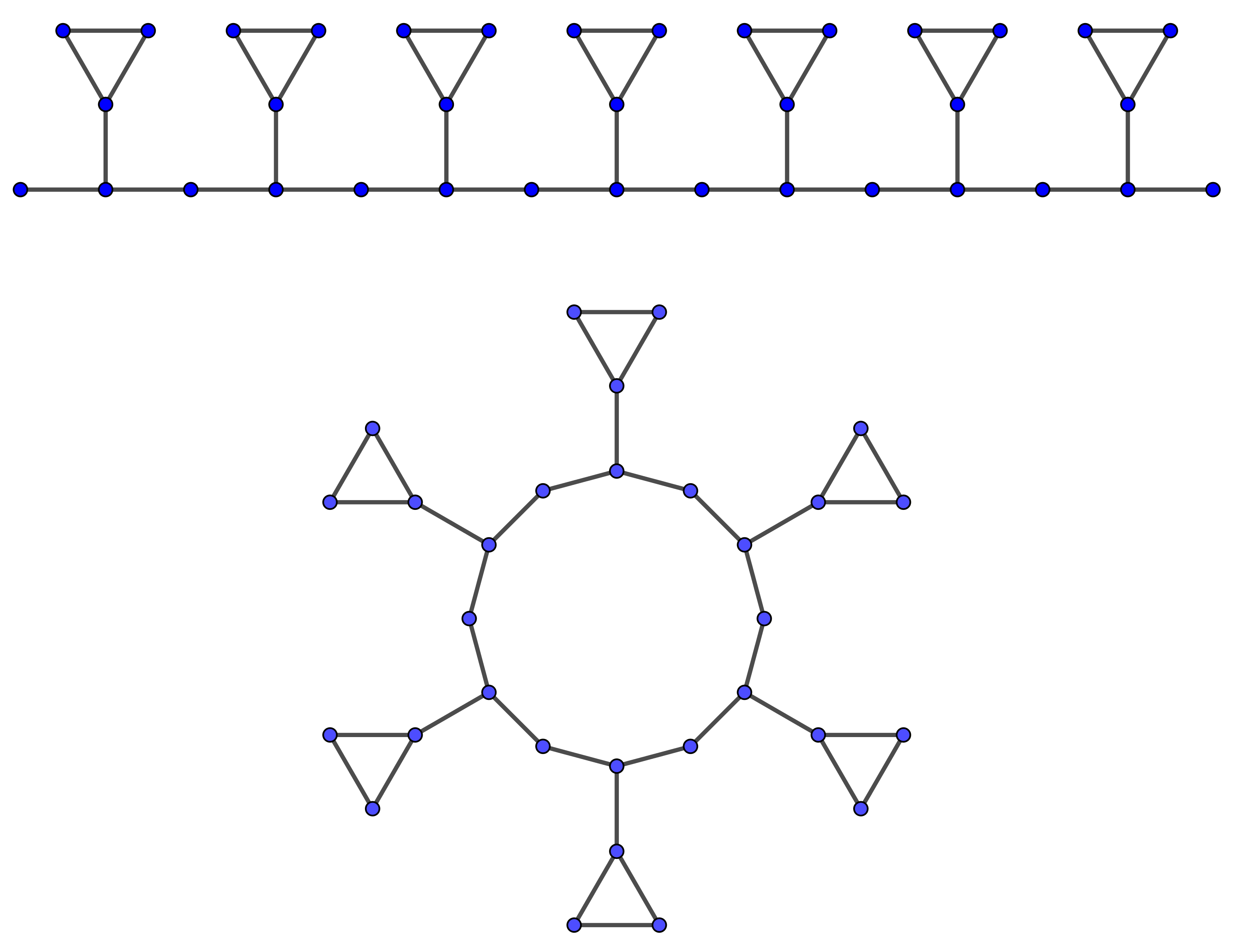}
\end{center}
\vspace{-6pt}
\caption{The wine glass path $W\!gp_7$ (top) and the wine glass cycle $W\!gc_6$ (bottom).}
\label{fig-wineglasses}
\end{figure}

We related the previous observation that 
the edge sets of most small counterexamples to Conjecture AAA 
may be partitioned into a tree-like disjoint union of triangles and 3-stars,
with the well-known facts that
the path has the largest energy among trees~\cite{gutman}
and that the cycle has either the largest or the second largest energy among unicyclic graphs~\cite{unicenergy1,unicenergy2},
to obtain (after several attempts) two families of graphs,
illustrated in Fig.~\ref{fig-wineglasses},
that will serve as infinite families of counterexamples to Conjecture AAA.
We aptly name these graphs as {\em wine glass paths} and {\em wine glass cycles}:
the wine glass path $W\!gp_k$ consists of $k$ ``wine glasses'' whose ``bases'' form a path,
while the wine glass cycle $W\!gc_k$ consists of $k$ wine glasses whose bases form a cycle.
$W\!gp_k$ has $5k+1$ vertices, $W\!gc_k$ has $5k$ vertices,
while both $W\!gp_k$ and $W\!gc_k$ have the matching number $\mu=2k$ and $\Delta=3$.

In order to prove that the graph sequences $(W\!gp_k)_{k \in \mathbb{N}}$ and $(W\!gc_k)_{k \in \mathbb{N} \setminus \{ 1 \} }$ contain infinitely many graphs $G$ such that $\mathcal{E}(G) \le 2 \mu(G)\sqrt{\Delta}$ does not hold, it is clearly sufficient to show that
\begin{align*}
    \lim_{k\to\infty} \frac{\mathcal{E}(W\!gp_k)}{\mu(W\!gp_k)} &> 2 \sqrt{3} \, ,\\
    \lim_{k\to\infty} \frac{\mathcal{E}(W\!gc_k)}{\mu(W\!gc_k)} &> 2 \sqrt{3} \, .
\end{align*}
We will do this in $3$ separate steps. First of all, we shall find the formulas for $\mathcal{E}(W\!gp_k)$ and $\mathcal{E}(W\!gc_k)$. Afterwards, we will use the obtained formulas in order to find a single expression for both $\lim\limits_{k\to\infty} \dfrac{\mathcal{E}(W\!gp_k)}{\mu(W\!gp_k)}$ and $\lim\limits_{k\to\infty} \dfrac{\mathcal{E}(W\!gc_k)}{\mu(W\!gc_k)}$, thus showing that these two values are equal. In the end, we will numerically compute the aforementioned expression up to a necessary precision and obtain a value greater than $2 \sqrt{3}$, which completes the entire proof.

Before we state the three main theorems which form the backbone of this section, we are going to need to prove an auxiliary lemma. As it will become clear soon, the function $\mathcal{F} \colon \mathbb{R} \setminus \{ -1, 2 \} \to \mathbb{R}$ defined by
$$\mathcal{F}(x) = x^2 - 3 - \dfrac{2}{(x+1)(x-2)}$$
is of great importance while calculating the energies of wine glass paths and wine glass cycles. We present the following lemma which proves some of its key properties:
\begin{lemma}\label{abcd_lemma}
    Let $y$ be an arbitrary real number such that $-2 \le y \le 2$. Then the following equation in $x \in \mathbb{R} \setminus \{-1, 2\}$
    $$ \mathcal{F}(x) = y$$
    has exactly $4$ distinct solutions $\alpha(y) < \beta(y) < \gamma(y) < \delta(y)$ such that
    \begin{align*}
        \alpha(y) &\in (-\infty, -1) \, ,\\
        \beta(y) &\in (-1, 0] \, ,\\
        \gamma(y) &\in (0, 2) \, ,\\
        \delta(y) &\in (2, +\infty) \, . 
    \end{align*}
    Moreover, we can interpret
    \begin{align*}
        \alpha &\colon [-2, 2] \to (-\infty, -1) \, ,\\
        \beta &\colon [-2, 2] \to (-1, 0] \, ,
    \end{align*}
    as real functions which are strictly decreasing and continuous, and
    \begin{align*}
        \gamma &\colon [-2, 2] \to (0, 2) \, ,\\
        \delta &\colon [-2, 2] \to (2, +\infty) \, ,
    \end{align*}
    as real functions which are strictly increasing and continuous.
\end{lemma}
\begin{proof}
    We can rewrite $\mathcal{F}(x)$ as
    $$\mathcal{F}(x) = x^2 - 3 + \dfrac{2}{3(x+1)} - \dfrac{2}{3(x-2)} \, .$$
    This function is clearly differentiable on its domain, which leads us to
    \begin{align*}
        \mathcal{F}'(x) &= 2x - \dfrac{2}{3(x+1)^2} + \dfrac{2}{3(x-2)^2}\\
        &= 2x + \dfrac{2(x+1)^2 - 2(x-2)^2}{3(x+1)^2(x-2)^2}\\
        &= 2x + \dfrac{2x^2 + 4x + 2 - 2x^2 + 8x - 8}{3(x+1)^2(x-2)^2}\\
        &= 2x + \dfrac{12x - 6}{3(x+1)^2(x-2)^2}\\
        &= 2x + \dfrac{4x - 2}{(x+1)^2(x-2)^2} \, .
    \end{align*}
    It becomes obvious that $\mathcal{F}'(x) > 0$ provided $x \ge \dfrac{1}{2}, x \neq 2$ and $\mathcal{F}'(x) < 0$ whenever $x \le 0,\ x \neq -1$. This means that the function $\mathcal{F}$ is strictly decreasing and continuous on $(-\infty, -1)$ and $(-1, 0]$ and strictly increasing and continuous on $(2, +\infty)$.
    
    Given the fact that
    \begin{align*}
        \lim_{x \to -\infty}\mathcal{F}(x) &= +\infty \, ,\\
        \lim_{x \to -1^{-}}\mathcal{F}(x) &= -\infty \, ,
    \end{align*}
    we conclude that for any $-2 \le y \le 2$ the equation $\mathcal{F}(x) = y$ must have a unique solution $\alpha(y)$ on $(-\infty, -1)$. Here, $\alpha \colon [-2, 2] \to (-\infty, -1)$ can be viewed as the restriction on $[-2, 2]$ of the inverse function of $\mathcal{F} \mid_{(-\infty, -1)} \colon (-\infty, -1) \to \mathbb{R}$, hence it must be both strictly decreasing and continuous. Similarly,
    \begin{align*}
        \lim_{x \to 2^+}\mathcal{F}(x) &= -\infty \, ,\\
        \lim_{x \to +\infty}\mathcal{F}(x) &= +\infty \, ,
    \end{align*}
    implies that for any $-2 \le y \le 2$, the equation $\mathcal{F}(x) = y$ must also have a unique solution $\delta(y)$ on $(2, +\infty)$. In this situation, $\delta \colon [-2, 2] \to (2, +\infty)$ can be interpreted as the restriction on $[-2, 2]$ of the inverse function of $\mathcal{F} \mid_{(2, +\infty)} \colon (2, +\infty) \to \mathbb{R}$, which means that it is a strictly increasing and continuous real function.

    By taking into consideration that $\mathcal{F}$ is strictly decreasing and continuous on $(-1, 0]$, together with
    \begin{align*}
        \lim_{x \to -1^{+}}\mathcal{F}(x) &= +\infty \, ,\\
        \mathcal{F}(0) &= -2 \, ,
    \end{align*}
    we further obtain that for all $-2 \le y \le 2$, the equation $\mathcal{F}(x) = y$ has a unique solution $\beta(y)$ on $(-1, 0]$. Due to the fact that $\beta \colon [-2, 2] \to (-1, 0]$ can be seen as the restriction on $[-2, 2]$ of the inverse function of $\mathcal{F} \mid_{(-1, 0]} \colon (-1, 0] \to [-2, +\infty)$, it is clear that this real function is strictly decreasing and continuous.

    In order to complete the proof, it is sufficient to prove that $\mathcal{F}(x) = y$ has exactly one solution on $(0, 2)$, for each $-2 \le y \le 2$, and then show the necessary properties of the corresponding real function. The function $\mathcal{F}'$ is differentiable on $(-1, 2)$, hence we get
    \begin{align*}
        \mathcal{F}''(x) &= 2 + \dfrac{4}{3(x+1)^3} - \dfrac{4}{3(x-2)^3} \, .
    \end{align*}
    For each $x \in (-1, 2)$, it can be directly seen that $\mathcal{F}''(x) > 0$, which means that $\mathcal{F}'$ is a strictly increasing function on $(-1, 2)$. We know that
    \begin{align*}
        \mathcal{F}'(0) &= 2 \cdot 0 - \dfrac{2}{3 \cdot 1^2} + \dfrac{2}{3 \cdot (-2)^2}\\
        &= -\dfrac{2}{3} + \dfrac{1}{6}\\
        &= -\dfrac{1}{2}
    \end{align*}
    and
    \begin{align*}
        \mathcal{F}'\left(\dfrac{1}{2}\right) &= 2 \cdot \dfrac{1}{2} - \dfrac{2}{3 \cdot \left( \frac{3}{2} \right)^2} + \dfrac{2}{3 \left( -\frac{3}{2} \right)^2}\\
        &= 1 - \dfrac{2}{\frac{27}{4}} + \dfrac{2}{\frac{27}{4}}\\
        &= 1 \, ,
    \end{align*}
    which implies that there must exist a real number $\xi \in \left(0, \dfrac{1}{2} \right)$ such that $\mathcal{F}'(\xi) = 0$. We then get $\mathcal{F}'(x) < 0$ for $x \in (-1, \xi)$ and $\mathcal{F}'(x) > 0$ for $x \in (\xi, 2)$, which shows that $\mathcal{F}$ is a strictly decreasing function on $(-1, \xi)$ and a strictly increasing function on $(\xi, 2)$.
    
    For $x \in (0, \xi]$, we know that $\mathcal{F}(x) < \mathcal{F}(0) = -2$, hence the equation $\mathcal{F}(x) = y$ cannot have any solutions for $-2 \le y \le 2$. The function $\mathcal{F}$ is strictly increasing and continuous on $(\xi, 2)$ and it is easy to see that
    \begin{align*}
        \mathcal{F}(\xi) &< -2 \, ,\\
        \lim_{x \to 2^{-}} \mathcal{F}(x) &= +\infty \, ,
    \end{align*}
    which leads us to conclude that for all $-2 \le y \le 2$ there must exist a unique solution $\gamma(y)$ to the equation $\mathcal{F}(x) = y$ on the interval $(\xi, 2)$. This means that $\gamma(y)$ is the unique solution on the interval $(0, 2)$ as well. We also notice that $\gamma \colon [-2, 2] \to (0, 2)$ can be regarded as the restriction on $[-2, 2]$ of the inverse function of $\mathcal{F}\mid_{(\xi, 2)} \colon (\xi, 2) \to (\mathcal{F}(\xi), +\infty)$, which shows that this function must be strictly increasing and continuous.
\end{proof}

In the remainder of the paper, we shall use $\alpha(y)$, $\beta(y)$, $\gamma(y)$ and $\delta(y)$ to denote the corresponding solutions to $\mathcal{F}(x) = y$, for any $-2 \le y \le 2$, as done so in Lemma \ref{abcd_lemma}. This allows us to state the three main theorems of this section, as follows:
\begin{theorem}\label{main_theorem_1}
    For any $k \ge 1$, the energy of the wine glass path $W\!gp_k$ is equal to
    \begin{equation}\label{main_formula_1}
        \mathcal{E}(W\!gp_k) = 2k - 2 \sum_{j=1}^{k} \alpha \left( 2 \cos\left( \frac{j \pi}{k+1}\right) \right) - 2 \sum_{j=1}^{k} \beta\left( 2 \cos\left( \frac{j \pi}{k+1}\right) \right) .
    \end{equation}
\end{theorem}
\begin{theorem}\label{main_theorem_2}
    For any $k \ge 2$, the energy of the wine glass cycle $W\!gc_k$ is equal to
    \begin{equation}\label{main_formula_2}
        \mathcal{E}(W\!gc_k) = 2k - 2 \sum_{j=0}^{k-1} \alpha \left( 2 \cos\left( \frac{2j \pi}{k}\right) \right) - 2 \sum_{j=0}^{k-1} \beta\left( 2 \cos\left( \frac{2j \pi}{k}\right) \right) .
    \end{equation}
\end{theorem}
\begin{theorem}\label{main_theorem_3}
    The graph sequences $(W\!gp_k)_{k \in \mathbb{N}}$ and $(W\!gc_k)_{k \in \mathbb{N} \setminus \{ 1 \} }$ satisfy
    $$\lim_{k\to\infty} \frac{\mathcal{E}(W\!gp_k)}{\mu(W\!gp_k)} = \lim_{k\to\infty} \frac{\mathcal{E}(W\!gc_k)}{\mu(W\!gc_k)} = L \, ,$$
    where
    \begin{equation}\label{main_formula_3}
        L = 1 - \alpha(-2) + \dfrac{1}{\pi} \int_{\alpha(2)}^{\alpha(-2)}\arccos\left(\dfrac{\mathcal{F}(x)}{2}\right)\diff x + \dfrac{1}{\pi} \int_{\beta(2)}^{0}\arccos\left(\dfrac{\mathcal{F}(x)}{2}\right)\diff x \, .
    \end{equation}
\end{theorem}

Before we begin writing out the proofs to these three theorems, we shall need one more auxiliary lemma.
\begin{lemma}\label{abcd_sum}
    For any $-2 \le y \le 2$, we have
    \begin{equation}\label{abcd_sum_eq}
        \alpha(y) + \beta(y) + \gamma(y) + \delta(y) = 1 \, .
    \end{equation}
\end{lemma}
\begin{proof}
    The $4$ distinct solutions $\alpha(y), \beta(y), \gamma(y), \delta(y)$ to $\mathcal{F}(x) = y$ must also be the solutions to the following equivalent equations in $x \in \mathbb{R} \setminus \{ -1, 2 \}$:
    \begin{alignat*}{2}
        && x^2 - 3 - \dfrac{2}{x^2-x-2} &= y\\
        \iff && (x^2-3)(x^2-x-2) - 2 &= y(x^2-x-2)\\
        \iff && x^4 - x^3 - 2x^2 - 3x^2 + 3x + 6 - 2 &= yx^2 - yx - 2y\\
        \iff && x^4 - x^3 - (y+5)x^2 + (y+3)x + (2y+4) &= 0 \, .
    \end{alignat*}
    Hence, we get that $\alpha(y), \beta(y), \gamma(y), \delta(y)$ actually represent the roots of the polynomial $x^4 - x^3 - (y+5)x^2 + (y+3)x + (2y+4)$. Eq.\ (\ref{abcd_sum_eq}) directly follows from Vieta's formulas.
\end{proof}

\bigskip\noindent
{\em Proof of Theorem \ref{main_theorem_1}}.\quad We will prove Eq.\ (\ref{main_formula_1}) by directly finding the set of all the eigenvalues of $W\!gp_k$, as well as the multiplicity of each eigenvalue.
Let $A$ denote the adjacency matrix of $W\!gp_k$. Given the fact that $A$ is necessarily a real symmetric matrix, all of its eigenvalues must be real. We conclude that in order to determine the spectrum of $W\!gp_k$, it is sufficient to solve the linear equation
\begin{equation}\label{wgp_eigen}
    Au = \lambda u
\end{equation}
in $u \in \mathbb{R}^{5k+1}$, parametrized by $\lambda \in \mathbb{R}$.
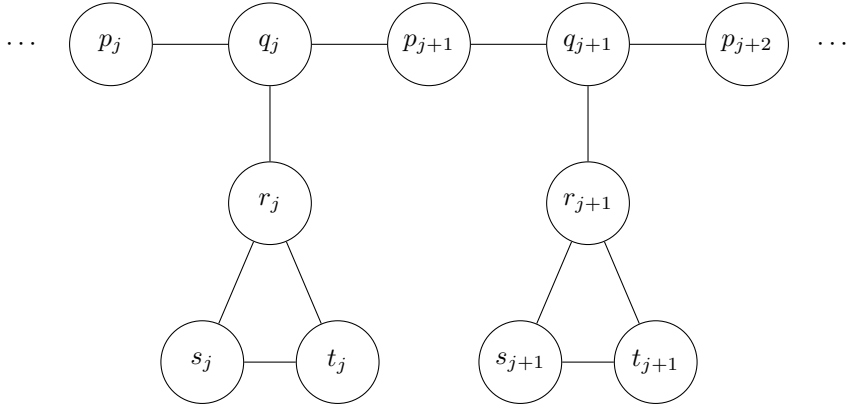
\begin{figure}[H]
    \centering
    \begin{tikzpicture}
        \node[state, minimum size=1.1cm] (1) {$p_j$};
        \node[state, minimum size=1.1cm] (2) [right =of 1] {$q_j$};
        \node[state, minimum size=1.1cm] (3) [below =of 2] {$r_j$};
        \node[state, minimum size=1.1cm] (4) [below =of 3, xshift=-0.9cm] {$s_j$};
        \node[state, minimum size=1.1cm] (5) [below =of 3, xshift=0.9cm] {$t_j$};
        \node[state, minimum size=1.1cm] (6) [right =of 2] {$p_{j+1}$};
        \node[state, minimum size=1.1cm] (7) [right =of 6] {$q_{j+1}$};
        \node[state, minimum size=1.1cm] (8) [below =of 7] {$r_{j+1}$};
        \node[state, minimum size=1.1cm] (9) [below =of 8, xshift=-0.9cm] {$s_{j+1}$};
        \node[state, minimum size=1.1cm] (10) [below =of 8, xshift=0.9cm] {$t_{j+1}$};
        \node[state, minimum size=1.1cm] (11) [right =of 7] {$p_{j+2}$};

        \node[minimum size=1.0cm] (12) [left =of 1, xshift=0.9cm] {$\cdots$};
        \node[minimum size=1.0cm] (13) [right =of 11, xshift=-0.9cm] {$\cdots$};

        \path (1) edge (2);
        \path (2) edge (3);
        \path (3) edge (4);
        \path (3) edge (5);
        \path (4) edge (5);
        \path (2) edge (6);
        \path (6) edge (7);
        \path (7) edge (8);
        \path (8) edge (9);
        \path (8) edge (10);
        \path (9) edge (10);
        \path (7) edge (11);
    \end{tikzpicture}
    \caption{The notation used to describe the elements of $u \in \mathbb{R}^{5k+1}$.}
    \label{wgp_notation}
\end{figure}
We will use the real values $p_j, j=\overline{0, k}$ and $q_j, r_j, s_j, t_j, j=\overline{0, k-1}$ to describe the $5k+1$ elements of $u$, in the manner done so in Figure \ref{wgp_notation}. Here it is important to notice that the vertices corresponding to the elements $p_0$ and $p_k$ are of degree $1$, and are only connected to the vertices corresponding to the elements $q_0$ and $q_{k-1}$, respectively.
Eq.\ (\ref{wgp_eigen}) promptly becomes
\begin{alignat}{2}
    \label{wgp_t}r_j + s_j &= \lambda t_j && \qquad (\forall j = \overline{0, k-1}) \, ,\\
    \label{wgp_s}r_j + t_j &= \lambda s_j && \qquad (\forall j = \overline{0, k-1}) \, ,\\
    \label{wgp_r}q_j + s_j + t_j &= \lambda r_j && \qquad (\forall j = \overline{0, k-1}) \, ,\\
    \label{wgp_q}p_j + r_j + p_{j+1} &= \lambda q_j && \qquad (\forall j = \overline{0, k-1}) \, ,\\
    \label{wgp_p}q_{j-1} + q_j &= \lambda p_j && \qquad (\forall j = \overline{1, k-1}) \, ,\\
    \label{wgp_p0}q_0 &= \lambda p_0 \, , && \\
    \label{wgp_pk}q_{k-1} &= \lambda p_k \, . &&
\end{alignat}
In order to solve the obtained equivalent system of equations, we will divide the problem into several cases.

\bigskip\noindent
\emph{Case }$\lambda = -1$.\quad From Eqs.\ (\ref{wgp_t}) and (\ref{wgp_s}) we get
$$r_j + s_j + t_j = 0$$
for each $0 \le j \le k-1$. Eq.\ (\ref{wgp_r}) gives
$$q_j + r_j + s_j + t_j = 0 \, ,$$
which immediately implies $q_j = 0$ for all $0 \le j \le k-1$. From Eqs.\ (\ref{wgp_p}), (\ref{wgp_p0}) and (\ref{wgp_pk}) we subsequently get that $p_j = 0$ for each $0 \le j \le k$. Finally, Eq.\ (\ref{wgp_q}) tells us that $r_j = 0$ holds as well, for all $0 \le j \le k-1$. Thus, we get that for $\lambda = -1$ each solution $u$ to Eq.\ (\ref{wgp_eigen}) satisfies
\begin{alignat}{2}
    \label{wgp_sol1}s_j + t_j &= 0 && \qquad (\forall j = \overline{0, k-1}) \, ,\\
    \label{wgp_sol2}r_j &= 0 && \qquad (\forall j = \overline{0, k-1}) \, ,\\
    \label{wgp_sol3}q_j &= 0 && \qquad (\forall j = \overline{0, k-1}) \, ,\\
    \label{wgp_sol4}p_j &= 0 && \qquad (\forall j = \overline{0, k}) \, .
\end{alignat}
It can easily be verified that each $u \in \mathbb{R}^{5k+1}$ satisfying Eqs.\ (\ref{wgp_sol1}), (\ref{wgp_sol2}), (\ref{wgp_sol3}) and (\ref{wgp_sol4}), is indeed a solution to Eq.\ (\ref{wgp_eigen}), provided $\lambda = -1$. We conclude that $-1$ must be an eigenvalue of $A$ and it is straightforward to see that its multiplicity equals $k$.

\bigskip\noindent
\emph{Case }$\lambda = 2$.\quad By subtracting Eq.\ (\ref{wgp_s}) from Eq.\ (\ref{wgp_t}) we get
\begin{alignat*}{2}
    && s_j - t_j &= \lambda (t_j - s_j)\\
    \implies \quad && (- 1 - \lambda)(t_j - s_j) &= 0\\
    \implies \quad && t_j - s_j &=0 \, ,
\end{alignat*}
which quickly implies that $s_j = t_j$ for each $0 \le j \le k-1$. Consequently, Eqs.\ (\ref{wgp_t}) and (\ref{wgp_s}) directly give $r_j = s_j = t_j$ for each $0 \le j \le k-1$. We now get from Eq.\ (\ref{wgp_r})
\begin{alignat*}{2}
    && q_j + s_j + t_j &= 2 r_j\\
    \implies \quad && q_j + 2 r_j &= 2 r_j\\
    \implies \quad && q_j &= 0 \, ,
\end{alignat*}
which means that $q_j = 0$ for all $0 \le j \le k-1$. Eqs.\ (\ref{wgp_p}), (\ref{wgp_p0}), (\ref{wgp_pk}) imply that $p_j = 0$ must also hold for each $0 \le j \le k$. However, Eq.\ (\ref{wgp_q}) now gives $r_j = 0$, which in turn leads us to $r_j = s_j = t_j = 0$ for all $0 \le j \le k-1$. Hence, we obtain $u = \mathbf{0}$. For $\lambda = 2$, Eq.\ (\ref{wgp_eigen}) only has the solution $\mathbf{0}$, so the conclusion is that $2$ is not an eigenvalue of $A$.

\bigskip\noindent
\emph{Case }$\lambda = 0$.\quad Eqs.\ (\ref{wgp_t}) and (\ref{wgp_s}) give $s_j = -r_j$ and $t_j = -r_j$ for each $0 \le j \le k-1$, respectively. Using Eq.\ (\ref{wgp_r}) we get
\begin{alignat*}{2}
    && q_j + s_j + t_j &= 0\\
    \implies \quad && q_j - r_j - r_j &= 0\\
    \implies \quad && q_j &= 2 r_j \, ,
\end{alignat*}
which means that $q_j = 2r_j$ for all $0 \le j \le k-1$. On the other hand, Eqs.\ (\ref{wgp_p0}) and (\ref{wgp_pk}) give us $q_0 = 0$ and $q_k = 0$. Together with Eq.\ (\ref{wgp_p}), this implies $q_j = 0$ for each $0 \le j \le k-1$. This subsequently gives $r_j = 0$, as well as $s_j = 0$ and $t_j = 0$. Taking into consideration Eq.\ (\ref{wgp_q}), we get that for $\lambda = 0$ each solution $u$ to Eq.\ (\ref{wgp_eigen}) must satisfy
\begin{alignat}{2}
    \label{wgp_sol5}p_j + p_{j+1} &= 0 && \qquad (\forall j = \overline{0, k-1}) \, ,\\
    \label{wgp_sol6}t_j &= 0 && \qquad (\forall j = \overline{0, k-1}) \, ,\\
    \label{wgp_sol7}s_j &= 0 && \qquad (\forall j = \overline{0, k-1}) \, ,\\
    \label{wgp_sol8}r_j &= 0 && \qquad (\forall j = \overline{0, k-1}) \, ,\\
    \label{wgp_sol9}q_j &= 0 && \qquad (\forall j = \overline{0, k-1}) \, .
\end{alignat}
It is trivial to check that each $u$ satisfying Eqs.\ (\ref{wgp_sol5}), (\ref{wgp_sol6}), (\ref{wgp_sol7}), (\ref{wgp_sol8}) and (\ref{wgp_sol9}) is indeed a solution to Eq.\ (\ref{wgp_eigen}), provided $\lambda = 0$. Since Eq.\ (\ref{wgp_sol5}) directly gives $p_j = (-1)^j p_0$ for all $0 \le j \le k$, we conclude that $0$ must be a simple eigenvalue of $A$.

\bigskip\noindent
\emph{Case }$\lambda \neq -1, 2, 0$.\quad If we subtract Eq.\ (\ref{wgp_s}) from Eq.\ (\ref{wgp_t}), we obtain
\begin{alignat*}{2}
    && s_j - t_j &= \lambda (t_j - s_j)\\
    \implies \quad && (- 1 - \lambda)(t_j - s_j) &= 0\\
    \implies \quad && t_j - s_j &=0\\
    \implies \quad && t_j &= s_j \, .
\end{alignat*}
Furthermore, Eq.\ (\ref{wgp_t}) implies $r_j = (\lambda - 1) t_j$ and Eq.\ (\ref{wgp_r}) gives us
\begin{alignat*}{2}
    && q_j + s_j + t_j &= \lambda r_j\\
    \implies \quad && q_j + 2t_j &= \lambda (\lambda - 1) t_j \\
    \implies \quad && q_j &= (\lambda^2 - \lambda - 2) t_j \\
    \implies \quad && q_j &= (\lambda+1)(\lambda-2) t_j \, .
\end{alignat*}
Hence, we conclude that $t_j = s_j = \dfrac{1}{(\lambda+1)(\lambda-2)} q_j$ and $r_j = \dfrac{\lambda-1}{(\lambda+1)(\lambda-2)} q_j$ for all $0 \le j \le k-1$. If we denote $q_{-1} = 0$ and $q_k = 0$, then from Eqs.\ (\ref{wgp_p}), (\ref{wgp_p0}) and (\ref{wgp_pk}) we also have
\begin{alignat*}{2}
    p_j &= \dfrac{1}{\lambda} q_{j-1} + \dfrac{1}{\lambda} q_j && \qquad (\forall j = \overline{0, k}) \, .
\end{alignat*}
Now we can use Eq.\ (\ref{wgp_q}) in order to obtain
\begin{alignat*}{2}
    && p_j + r_j + p_{j+1} &= \lambda q_j\\
    \implies \quad && \dfrac{q_{j-1}}{\lambda} + \dfrac{q_j}{\lambda} + \dfrac{(\lambda-1) q_j}{(\lambda+1)(\lambda-2)} + \dfrac{q_j}{\lambda} + \dfrac{q_{j+1}}{\lambda} &= \lambda q_j\\
    \implies \quad && q_{j-1} + q_j + \dfrac{\lambda(\lambda-1) q_j}{(\lambda+1)(\lambda-2)}  + q_j + q_{j+1} &= \lambda^2 q_j\\
    \implies \quad && q_{j-1} + q_{j+1} &= \left( \lambda^2 - 2 - \dfrac{\lambda(\lambda-1)}{(\lambda+1)(\lambda-2)} \right) q_j\\
    \implies \quad && q_{j-1} + q_{j+1} &= \left( \lambda^2 - 3 - \dfrac{2}{(\lambda + 1)(\lambda + 2)} \right) q_j\\
    \implies \quad && q_{j-1} + q_{j+1} &= \mathcal{F}(\lambda) \ q_j
\end{alignat*}
for all $0 \le j \le k-1$. Thus, we get that whenever $\lambda \neq -1, 2, 0$, each solution $u$ to Eq.\ (\ref{wgp_eigen}) satisfies
\begin{alignat}{2}
    \label{wgp_sol10} t_j &= \dfrac{1}{(\lambda+1)(\lambda-2)} q_j && \qquad (\forall j = \overline{0, k-1}) \, ,\\
    \label{wgp_sol11} s_j &= \dfrac{1}{(\lambda+1)(\lambda-2)} q_j && \qquad (\forall j = \overline{0, k-1}) \, ,\\
    \label{wgp_sol12} r_j &= \dfrac{\lambda-1}{(\lambda+1)(\lambda-2)} q_j && \qquad (\forall j = \overline{0, k-1}) \, ,\\
    \label{wgp_sol13}p_j &= \dfrac{1}{\lambda} q_{j-1} + \dfrac{1}{\lambda} q_j && \qquad (\forall j = \overline{0, k}) \, ,\\
    \label{wgp_sol14}q_{j-1} + q_{j+1} &= \mathcal{F}(\lambda) \ q_j && \qquad (\forall j = \overline{0, k-1}) \, ,
\end{alignat}
where $q_{-1} = q_k = 0$, as noted earlier. The converse can also directly be shown, i.e.\ each $u$ such that Eqs.\ (\ref{wgp_sol10}), (\ref{wgp_sol11}), (\ref{wgp_sol12}), (\ref{wgp_sol13}) and (\ref{wgp_sol14}) all hold, must be a solution to Eq.\ (\ref{wgp_eigen}), provided $\lambda \neq -1, 2, 0$.

If we denote $w = \begin{bmatrix} q_0 & q_1 & \cdots & q_{k-1} \end{bmatrix}^T$, we then see that Eq.\ (\ref{wgp_sol14}) becomes equivalent to
\begin{align}\label{wgp_sol15}
    Bw &= \mathcal{F}(\lambda) \ w ,
\end{align}
where
\begin{align*}
    B = \begin{bmatrix}
    0 & 1 & 0 & 0 & \cdots & 0 & 0\\
    1 & 0 & 1 & 0 & \cdots & 0 & 0\\
    0 & 1 & 0 & 1 & \cdots & 0 & 0\\
    0 & 0 & 1 & 0 & \cdots & 0 & 0\\
    \vdots & \vdots & \vdots & \vdots & \ddots & \vdots & \vdots\\
    0 & 0 & 0 & 0 & \cdots & 0 & 1\\
    0 & 0 & 0 & 0 & \cdots & 1 & 0
\end{bmatrix} \in \mathbb{R}^{k \times k} \, .
\end{align*}
Here, we can notice that $B$ actually represents the adjacency matrix of a path graph on $k$ vertices. It is known (see, for example, \cite[pp.\ 18]{spectra_of_graphs}) that the adjacency matrix of a path graph on $k$ vertices must have $k$ distinct simple eigenvalues $2\cos\left(\dfrac{1}{k+1}\pi\right), 2\cos\left(\dfrac{2}{k+1}\pi\right), \ldots, 2\cos\left(\dfrac{k}{k+1}\pi\right)$.

Suppose that $\lambda \neq -1, 2, 0$ is such that $\mathcal{F}(\lambda) = 2\cos\left(\dfrac{j}{k+1}\pi\right)$ for some $1 \le j \le k$. We then get that there exists a $w \neq \mathbf{0}$ such that Eq.\ (\ref{wgp_sol15}) holds. This basically means that there exist some $q_0, q_1, \ldots, q_{k-1}$, at least one of which is nonzero, which satisfy Eq.\ (\ref{wgp_sol14}). By using Eqs.\ (\ref{wgp_sol10}), (\ref{wgp_sol11}), (\ref{wgp_sol12}) and (\ref{wgp_sol13}) in order to construct the other elements of $u$, we conclude that there must exist some $u \neq \mathbf{0}$ such that $u$ is a solution to Eq.\ (\ref{wgp_eigen}). Thus, $\lambda$ must be an eigenvalue of $A$.

All the values $2\cos\left(\dfrac{j}{k+1}\pi\right)$ for $1 \le j \le k$ lie within the interval $(-2, 2)$. By using Lemma \ref{abcd_lemma}, we see that there exist $4$ distinct real values $\alpha\left(2\cos\left(\dfrac{j}{k+1}\pi\right)\right), \linebreak \beta\left( 2\cos\left(\dfrac{j}{k+1}\pi\right)\right), \gamma\left(2\cos\left(\dfrac{j}{k+1}\pi\right)\right), \delta\left(2\cos\left(\dfrac{j}{k+1}\pi\right)\right)$ which represent \linebreak the solutions to the equation $\mathcal{F}(x) = 2\cos\left(\dfrac{j}{k+1}\pi\right)$ in $x \in \mathbb{R}$, for each $1 \le j \le k$. These $4k$ real numbers are clearly all distinct and none of them are equal to $-1$ or $2$, due to Lemma \ref{abcd_lemma}. It can easily be checked that none of them are equal to $0$ either, since $\mathcal{F}(0) = -2$ and all of the $2\cos\left(\dfrac{j}{k+1}\pi\right)$ are greater than $-2$. We conclude that all of the elements of the set
\begin{align*}
    Z = \{ \alpha(y), \beta(y), \gamma(y), \delta(y) \colon y = 2\cos\left(\dfrac{j}{k+1}\pi\right), j=\overline{1, k} \}
\end{align*}
must be an eigenvalue of $A$. Bearing in mind that the matrix $A$ is of order $5k+1$ and we already know that it has an eigenvalue $-1$ of multiplicity $k$ and a simple eigenvalue $0$, we see that all of its eigenvalues from the set $Z$ must be simple eigenvalues. Also, there can be more no eigenvalues other than the ones we have already found, which means that we have fully determined the spectrum of $A$. As a direct consequence, we obtain the formula
\begin{align*}
    \mathcal{E}(W\!gp_k) = k\lvert -1 \rvert &+ \sum_{j=1}^k \left| \alpha\left(2\cos\left(\dfrac{j}{k+1}\pi\right)\right) \right| + \sum_{j=1}^k \left| \beta\left(2\cos\left(\dfrac{j}{k+1}\pi\right)\right) \right|\\
    &+ \sum_{j=1}^k \left| \gamma\left(2\cos\left(\dfrac{j}{k+1}\pi\right)\right) \right| + \sum_{j=1}^k \left| \delta\left(2\cos\left(\dfrac{j}{k+1}\pi\right)\right) \right| .
\end{align*}
Since $\alpha(y) < \beta(y) \le 0$ and $0 < \gamma(y) < \delta(y)$ for all $-2 \le y \le 2$, we get
\begin{align*}
    \mathcal{E}(W\!gp_k) = k &- \sum_{j=1}^k \alpha\left(2\cos\left(\dfrac{j}{k+1}\pi\right)\right) - \sum_{j=1}^k \beta\left(2\cos\left(\dfrac{j}{k+1}\pi\right)\right)\\
    &+ \sum_{j=1}^k \gamma\left(2\cos\left(\dfrac{j}{k+1}\pi\right)\right) + \sum_{j=1}^k  \delta\left(2\cos\left(\dfrac{j}{k+1}\pi\right)\right) .
\end{align*}
Finally, Lemma \ref{abcd_sum} gives
\begin{alignat*}{2}
    && \alpha(y) + \beta(y) + \gamma(y) + \delta(y) &= 1\\
    \implies \quad && -\alpha(y)-\beta(y)+\gamma(y)+\delta(y) &= 1 - 2\alpha(y) - 2\beta(y) \, ,
\end{alignat*}
for each $-2 \le y \le 2$, which directly leads us to
\begin{align*}
    \mathcal{E}(W\!gp_k) &= k + \sum_{j=1}^k 1 - \sum_{j=1}^k 2 \alpha\left(2\cos\left(\dfrac{j}{k+1}\pi\right)\right) - \sum_{j=1}^k 2 \beta\left(2\cos\left(\dfrac{j}{k+1}\pi\right)\right)\\
    &= 2k - 2 \sum_{j=1}^k \alpha\left(2\cos\left(\dfrac{j}{k+1}\pi\right)\right) - 2 \sum_{j=1}^k \beta\left(2\cos\left(\dfrac{j}{k+1}\pi\right)\right) .\qed
\end{align*}

\bigskip\noindent
{\em Proof of Theorem \ref{main_theorem_2}}.\quad This proof is also based on directly finding the set of all the eigenvalues of the corresponding graph $W\!gc_k$, where $k \ge 2$, together with all of their multiplicities. Let $A$ be the adjacency matrix of $W\!gc_k$. Due to the fact that $A$ is a real symmetric matrix, all of its eigenvalues must be real. This means that we can fully determine the spectrum of $A$ by solving the linear equation
\begin{align}\label{wgc_eigen}
    Au &= \lambda u
\end{align}
in $u \in \mathbb{R}^{5k}$, parametrized by $\lambda \in \mathbb{R}$.
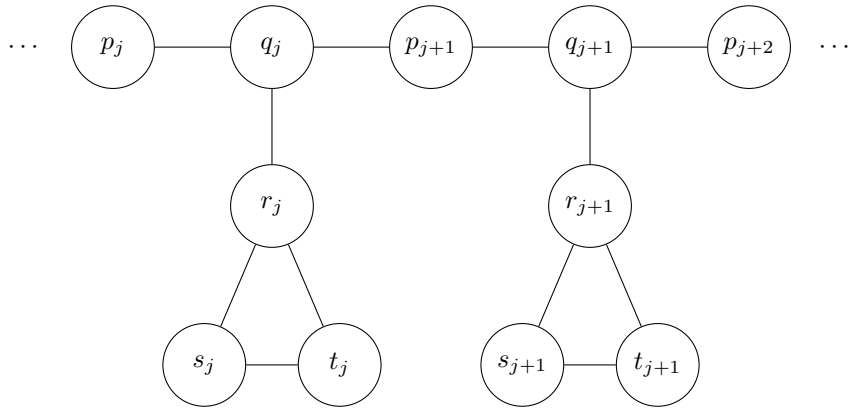
\begin{figure}[H]
    \centering
    \begin{tikzpicture}
        \node[state, minimum size=1.1cm] (1) {$p_j$};
        \node[state, minimum size=1.1cm] (2) [right =of 1] {$q_j$};
        \node[state, minimum size=1.1cm] (3) [below =of 2] {$r_j$};
        \node[state, minimum size=1.1cm] (4) [below =of 3, xshift=-0.9cm] {$s_j$};
        \node[state, minimum size=1.1cm] (5) [below =of 3, xshift=0.9cm] {$t_j$};
        \node[state, minimum size=1.1cm] (6) [right =of 2] {$p_{j+1}$};
        \node[state, minimum size=1.1cm] (7) [right =of 6] {$q_{j+1}$};
        \node[state, minimum size=1.1cm] (8) [below =of 7] {$r_{j+1}$};
        \node[state, minimum size=1.1cm] (9) [below =of 8, xshift=-0.9cm] {$s_{j+1}$};
        \node[state, minimum size=1.1cm] (10) [below =of 8, xshift=0.9cm] {$t_{j+1}$};
        \node[state, minimum size=1.1cm] (11) [right =of 7] {$p_{j+2}$};

        \node[minimum size=1.0cm] (12) [left =of 1, xshift=0.9cm] {$\cdots$};
        \node[minimum size=1.0cm] (13) [right =of 11, xshift=-0.9cm] {$\cdots$};

        \path (1) edge (2);
        \path (2) edge (3);
        \path (3) edge (4);
        \path (3) edge (5);
        \path (4) edge (5);
        \path (2) edge (6);
        \path (6) edge (7);
        \path (7) edge (8);
        \path (8) edge (9);
        \path (8) edge (10);
        \path (9) edge (10);
        \path (7) edge (11);
    \end{tikzpicture}
    \caption{The notation used to describe the elements of $u \in \mathbb{R}^{5k}$.}
    \label{wgc_notation}
\end{figure}
We will use the real values $p_j, q_j, r_j, s_j, t_j, j=\overline{0, k-1}$ to describe the $5k$ elements of $u$, according to Figure \ref{wgc_notation}. Here, the vertices corresponding to the elements $q_{k-1}$ and $p_0$ are connected as well, due to the cyclic nature of the graph. Eq.\ (\ref{wgc_eigen}) now becomes equivalent to
\begin{alignat}{2}
    \label{wgc_t}r_j + s_j &= \lambda t_j && \qquad (\forall j = \overline{0, k-1}) \, ,\\
    \label{wgc_s}r_j + t_j &= \lambda s_j && \qquad (\forall j = \overline{0, k-1}) \, ,\\
    \label{wgc_r}q_j + s_j + t_j &= \lambda r_j && \qquad (\forall j = \overline{0, k-1}) \, ,\\
    \label{wgc_q}p_j + r_j + p_{j+1} &= \lambda q_j && \qquad (\forall j = \overline{0, k-1}) \, ,\\
    \label{wgc_p}q_{j-1} + q_j &= \lambda p_j && \qquad (\forall j = \overline{0, k-1}) \, ,
\end{alignat}
where we define $p_k = p_0$ and $q_{-1} = q_{k-1}$. We will solve the given parametrized system of equations by dividing the problem into several cases, as done so in the proof of Theorem \ref{main_theorem_1}.

\bigskip\noindent
\emph{Case }$\lambda = -1$.\quad Eqs.\ (\ref{wgc_t}) and (\ref{wgc_s}) give us
$$r_j + s_j + t_j = 0$$
for each $0 \le j \le k-1$, while Eq.\ (\ref{wgc_r}) implies
$$q_j + r_j + s_j + t_j = 0$$
for all $0 \le j \le k-1$. This means that all the $q_j$ must be equal to zero, hence all the $p_j$ must be equal to zero as well, due to Eq. (\ref{wgc_p}). We now get from Eq.\ (\ref{wgc_q}) that $r_j = 0$ for each $0 \le j \le k-1$, which leads us to conclude that for $\lambda = -1$ each solution $u$ to Eq.\ (\ref{wgc_eigen}) must satisfy
\begin{alignat}{2}
    \label{wgc_sol1}s_j + t_j &= 0 && \qquad (\forall j = \overline{0, k-1}) \, ,\\
    \label{wgc_sol2}r_j &= 0 && \qquad (\forall j = \overline{0, k-1}) \, ,\\
    \label{wgc_sol3}q_j &= 0 && \qquad (\forall j = \overline{0, k-1}) \, ,\\
    \label{wgc_sol4}p_j &= 0 && \qquad (\forall j = \overline{0, k-1}) \, .
\end{alignat}
It can be directly shown that the converse is also true, i.e.\ each $u \in \mathbb{R}^{5k}$ which satisfies (\ref{wgc_sol1}), (\ref{wgc_sol2}), (\ref{wgc_sol3}) and (\ref{wgc_sol4}) must also be a solution to Eq.\ (\ref{wgc_eigen}), provided $\lambda = -1$. We obtain that $-1$ must be an eigenvalue of $A$ and it is easy to see that its multiplicity equals $k$.

\bigskip\noindent
\emph{Case }$\lambda = 2$.\quad If we subtract Eq.\ (\ref{wgc_s}) from Eq.\ (\ref{wgc_t}), we have
\begin{alignat*}{2}
    && s_j - t_j &= \lambda(t_j - s_j)\\
    \implies \quad && (-\lambda-1)(t_j - s_j) &= 0\\
    \implies \quad && t_j - s_j &= 0 \, .
\end{alignat*}
Hence, $s_j = t_j$ for all $0 \le j \le k-1$. Eqs.\ (\ref{wgc_t}) and (\ref{wgc_s}) now give us $r_j = s_j = t_j$ for all $0 \le j \le k-1$. By using Eq.\ (\ref{wgc_r}), we further obtain
\begin{alignat*}{2}
    && q_j + s_j + t_j &= 2r_j\\
    \implies \quad && q_j + 2r_j &= 2r_j\\
    \implies \quad && q_j &= 0 \, .
\end{alignat*}
We get that all the $q_j$ must be equal to zero, which directly implies that all the $p_j$ are equal to zero as well, due to Eq.\ (\ref{wgc_p}). However, Eq.\ (\ref{wgc_q}) now gives us $r_j = 0$ for each $0 \le j \le k-1$. This implies $s_j = t_j = 0$ as well, hence $u = \mathbf{0}$. We conclude that for $\lambda = 2$, the only solution to Eq.\ (\ref{wgc_eigen}) is $u = \mathbf{0}$, which means that $2$ is not an eigenvalue of $A$.

\bigskip\noindent
\emph{Case }$\lambda = 0$.\quad Eqs.\ (\ref{wgc_t}) and (\ref{wgc_s}) directly give $s_j = t_j = -r_j$ for all $0 \le j \le k-1$. We then get from Eq.\ (\ref{wgc_r})
\begin{alignat*}{2}
    && q_j + s_j + t_j &= 0\\
    \implies \quad && q_j - 2 r_j &= 0\\
    \implies \quad && q_j &= 2r_j ,
\end{alignat*}
hence $q_j = 2r_j$ holds for each $0 \le j \le k-1$. Also, Eq.\ (\ref{wgc_p}) implies $q_j = (-1)^j q_0$ for all $0 \le j \le k-1$. The further logical deductions depend on the parity of $k$, so we break down the problem into two subcases.

\noindent
\emph{Subcase }$k$ is odd.\quad We have
\begin{alignat*}{2}
    && q_{k-1} + q_0 &= 0\\
    \implies \quad && ((-1)^{k-1} + 1) q_0 &= 0\\
    \implies \quad && 2q_0 &= 0\\
    \implies \quad && q_0 &= 0 ,
\end{alignat*}
which means that all the $q_j$ must be equal to zero. Thus, we obtain $q_j = r_j = s_j = t_j = 0$ for all $0 \le j \le k-1$. Eq.\ (\ref{wgc_q}) now becomes
\begin{alignat*}{2}
    p_j + p_{j+1} &= 0 && \qquad (\forall j = \overline{0, k-1}) \, .
\end{alignat*}
This means that $p_j = (-1)^j p_0$ for each $0 \le j \le k-1$, which implies
\begin{alignat*}{2}
    && p_{k-1} + p_0 &= 0\\
    \implies \quad && ((-1)^{k-1} + 1) p_0 &= 0\\
    \implies \quad && 2p_0 &= 0\\
    \implies \quad && p_0 &= 0 .
\end{alignat*}
We have proved that $u = \mathbf{0}$ is the unique solution to Eq.\ (\ref{wgc_eigen}), provided $\lambda = 0$ and $k$ is odd. This means that $0$ is not an eigenvalue of $A$ when $k$ is odd.

\noindent
\emph{Subcase }$k$ is even.\quad From $q_j = (-1)^j q_0$ and $q_j = 2r_j$, it is clear that $r_j = (-1)^j r_0$ for all $0 \le j \le k-1$. By summing Eq.\ (\ref{wgc_q}) for $j = 0, 2, 4, \ldots, k-2$, we get
\begin{alignat*}{2}
    && \sum_{j=0}^{k-1}p_j + \sum_{j=0}^{(k-2)/2}r_{2j} &= 0\\
    \implies \quad && \sum_{j=0}^{k-1}p_j + \frac{k}{2}r_0 &= 0\\
    \implies \quad && \sum_{j=0}^{k-1}p_j &= -\frac{k}{2}r_0 \, .
\end{alignat*}
By summing the same equation for $j = 1, 3, 5, \ldots, k-1$, we obtain
\begin{alignat*}{2}
    && \sum_{j=0}^{k-1}p_j + \sum_{j=0}^{(k-2)/2}r_{2j+1} &= 0\\
    \implies \quad && \sum_{j=0}^{k-1}p_j + \frac{k}{2}r_1 &= 0\\
    \implies \quad && \sum_{j=0}^{k-1}p_j &= -\frac{k}{2}r_1 \, ,
\end{alignat*}
which helps us conclude that $r_0 = r_1$. Since $r_1 = -r_0$, this means that $r_0 = 0$, which promptly implies $q_j = r_j = s_j = t_j = 0$ for all $0 \le j \le k-1$. Eq.\ (\ref{wgc_q}) now becomes
\begin{alignat*}{2}
    p_j + p_{j+1} &= 0 && \qquad (\forall j = \overline{0, k-1}) \, ,
\end{alignat*}
hence $p_j = (-1)^j p_0$ for each $0 \le j \le k-1$. Thus, we get that when $\lambda = 0$ and $k$ is even, each solution to Eq.\ (\ref{wgc_eigen}) satisfies
\begin{alignat}{2}
    \label{wgc_sol5}t_j &= 0 && \qquad (\forall j = \overline{0, k-1}) \, ,\\
    \label{wgc_sol6}s_j &= 0 && \qquad (\forall j = \overline{0, k-1}) \, ,\\
    \label{wgc_sol7}r_j &= 0 && \qquad (\forall j = \overline{0, k-1}) \, ,\\
    \label{wgc_sol8}q_j &= 0 && \qquad (\forall j = \overline{0, k-1}) \, ,\\
    \label{wgc_sol9}p_j &= (-1)^j p_0 && \qquad (\forall j = \overline{0, k-1}) \, .
\end{alignat}
It can directly be shown that if $p_0$ is chosen as an arbitrary real value, and all the other elements of $u$ are determined according to Eqs.\ (\ref{wgc_sol5}), (\ref{wgc_sol6}), (\ref{wgc_sol7}), (\ref{wgc_sol8}) and (\ref{wgc_sol9}), then the obtained $u$ must be a solution to Eq.\ (\ref{wgc_eigen}), provided $\lambda = 0$ and $k$ is even. We conclude that $0$ is a simple eigenvalue of $A$ when $k$ is even.

\bigskip\noindent
\emph{Case }$\lambda \neq -1, 2, 0$.\quad If we subtract Eq.\ (\ref{wgc_s}) from Eq.\ (\ref{wgc_t}), we get
\begin{alignat*}{2}
    && s_j - t_j &= \lambda (t_j - s_j)\\
    \implies \quad && (- 1 - \lambda)(t_j - s_j) &= 0\\
    \implies \quad && t_j - s_j &=0\\
    \implies \quad && t_j &= s_j \, .
\end{alignat*}
At the same time, Eq.\ (\ref{wgc_t}) gives $r_j = (\lambda - 1)t_j$ and from Eq.\ (\ref{wgc_r}) we obtain
\begin{alignat*}{2}
    && q_j + s_j + t_j &= \lambda r_j\\
    \implies \quad && q_j + 2t_j &= \lambda (\lambda - 1) t_j \\
    \implies \quad && q_j &= (\lambda^2 - \lambda - 2) t_j \\
    \implies \quad && q_j &= (\lambda+1)(\lambda-2) t_j \, .
\end{alignat*}
This means that $t_j = s_j = \dfrac{1}{(\lambda+1)(\lambda-2)} q_j$ and $r_j = \dfrac{\lambda-1}{(\lambda+1)(\lambda-2)} q_j$ for all $0 \le j \le k-1$. On the other hand, Eq.\ (\ref{wgc_p}) implies
\begin{alignat*}{2}
    p_j &= \dfrac{1}{\lambda} q_{j-1} + \dfrac{1}{\lambda} q_j && \qquad (\forall j = \overline{0, k-1}) \, ,
\end{alignat*}
which allows us to denote $q_k = q_0$ and use Eq.\ (\ref{wgc_q}) in order to compute
\begin{alignat*}{2}
    && p_j + r_j + p_{j+1} &= \lambda q_j\\
    \implies \quad && \dfrac{q_{j-1}}{\lambda} + \dfrac{q_j}{\lambda} + \dfrac{(\lambda-1) q_j}{(\lambda+1)(\lambda-2)} + \dfrac{q_j}{\lambda} + \dfrac{q_{j+1}}{\lambda} &= \lambda q_j\\
    \implies \quad && q_{j-1} + q_j + \dfrac{\lambda(\lambda-1) q_j}{(\lambda+1)(\lambda-2)}  + q_j + q_{j+1} &= \lambda^2 q_j\\
    \implies \quad && q_{j-1} + q_{j+1} &= \left( \lambda^2 - 2 - \dfrac{\lambda(\lambda-1)}{(\lambda+1)(\lambda-2)} \right) q_j\\
    \implies \quad && q_{j-1} + q_{j+1} &= \left( \lambda^2 - 3 - \dfrac{2}{(\lambda + 1)(\lambda + 2)} \right) q_j\\
    \implies \quad && q_{j-1} + q_{j+1} &= \mathcal{F}(\lambda) \ q_j \, .
\end{alignat*}
We conclude that, provided $\lambda \neq -1, 2, 0$, every solution $u$ to Eq.\ (\ref{wgc_eigen}) must satisfy
\begin{alignat}{2}
    \label{wgc_sol10} t_j &= \dfrac{1}{(\lambda+1)(\lambda-2)} q_j && \qquad (\forall j = \overline{0, k-1}) \, ,\\
    \label{wgc_sol11} s_j &= \dfrac{1}{(\lambda+1)(\lambda-2)} q_j && \qquad (\forall j = \overline{0, k-1}) \, ,\\
    \label{wgc_sol12} r_j &= \dfrac{\lambda-1}{(\lambda+1)(\lambda-2)} q_j && \qquad (\forall j = \overline{0, k-1}) \, ,\\
    \label{wgc_sol13}p_j &= \dfrac{1}{\lambda} q_{j-1} + \dfrac{1}{\lambda} q_j && \qquad (\forall j = \overline{0, k-1}) \, ,\\
    \label{wgc_sol14}q_{j-1} + q_{j+1} &= \mathcal{F}(\lambda) \ q_j && \qquad (\forall j = \overline{0, k-1}) \, .
\end{alignat}
The converse is also directly shown, i.e.\ every $u$ such that Eqs.\ (\ref{wgc_sol10}), (\ref{wgc_sol11}), (\ref{wgc_sol12}), (\ref{wgc_sol13}) and (\ref{wgc_sol14}) all hold, must be a solution to Eq.\ (\ref{wgc_eigen}), in case $\lambda \neq -1, 2, 0$.

Let us denote $w = \begin{bmatrix} q_0 & q_1 & \cdots & q_{k-1} \end{bmatrix}^T$. It is clear that Eq.\ (\ref{wgc_sol14}) is equivalent to
\begin{align}\label{wgc_sol15}
    Bw &= \mathcal{F}(\lambda) \ w ,
\end{align}
where
\begin{align*}
    B = \begin{bmatrix}
    0 & 1 & 0 & 0 & \cdots & 0 & 1\\
    1 & 0 & 1 & 0 & \cdots & 0 & 0\\
    0 & 1 & 0 & 1 & \cdots & 0 & 0\\
    0 & 0 & 1 & 0 & \cdots & 0 & 0\\
    \vdots & \vdots & \vdots & \vdots & \ddots & \vdots & \vdots\\
    0 & 0 & 0 & 0 & \cdots & 0 & 1\\
    1 & 0 & 0 & 0 & \cdots & 1 & 0
\end{bmatrix} \in \mathbb{R}^{k \times k} \, .
\end{align*}
It is important to notice that $B$ represents the adjacency matrix of a cycle graph on $k$ vertices. It is known (see, for example, \cite[pp.\ 18]{spectra_of_graphs}) that the spectrum of the adjacency matrix of a cycle graph on $k$ vertices must be composed of the real numbers $2\cos\left(\dfrac{2j \pi}{k} \right)$ for $j = \overline{0, k-1}$.

Suppose that $\lambda \neq -1, 2, 0$ is such that $\mathcal{F}(\lambda)$ is an eigenvalue of $B$ with a multiplicity of $g$. In this case there exists a set of linearly independent vectors $\{w_1, w_2, \ldots, w_g \}$ which are all solutions to Eq.\ (\ref{wgc_sol15}) for the selected value of $\lambda$. Furthermore, for each vector $w_j$, $1 \le j \le g$, we can construct the corresponding vector $u_j$ by using Eqs.\ (\ref{wgc_sol10}), (\ref{wgc_sol11}), (\ref{wgc_sol12}) and (\ref{wgc_sol13}). It is easy to see that each $u_j$ will be a solution to Eq.\ (\ref{wgc_eigen}), for $1 \le j \le g$. An important observation is that the matrix
\begin{align*}
    W = \left[
    \begin{array}{c|c|c|c}
        w_1 & w_2 & \cdots & w_g
    \end{array}
    \right]
\end{align*}
will be a submatrix of the matrix
\begin{align*}
    U = \left[
    \begin{array}{c|c|c|c}
        u_1 & u_2 & \cdots & u_g
    \end{array}
    \right] .
\end{align*}
$W$ has $g$ linearly independent columns, hence its rank must be equal to $g$. This means that $U$ cannot have a rank lower than $g$. However, since $U$ also has exactly $g$ columns, we conclude that its rank is equal to $g$ as well. We obtain that the columns of $U$ must be linearly independent as well, which tells us that $\lambda$ must be an eigenvalue of $A$ with a multiplicity not smaller than $g$.

The further logical deductions slightly differ depending on the parity of $k$, so we divide the problem into two subcases.

\noindent
\emph{Subcase }$k$ is odd.\quad Let $k = 2k_1 + 1$. The matrix $B$ must have a simple eigenvalue $2$, as well as $k_1$ additional eigenvalues $2\cos\left(\dfrac{2\pi}{k} \right), 2\cos\left(\dfrac{4\pi}{k} \right), \ldots, 2\cos\left(\dfrac{2k_1 \pi}{k} \right)$, each of which has a multiplicity of $2$.

All of the values $2\cos\left(\dfrac{2j \pi}{k} \right)$ for $1 \le j \le k_1$ lie within the interval $(-2, 2)$. According to Lemma \ref{abcd_lemma}, there exist $4$ distinct real numbers $\alpha\left(2\cos\left(\dfrac{2j \pi}{k}\right)\right), \linebreak \beta\left(2\cos\left(\dfrac{2j \pi}{k}\right)\right), \gamma\left(2\cos\left(\dfrac{2j \pi}{k}\right)\right), \delta\left(2\cos\left(\dfrac{2j \pi}{k}\right)\right)$ which represent the solutions to the equation $\mathcal{F}(x) = 2\cos\left(\dfrac{2j \pi}{k}\right)$ in $x \in \mathbb{R}$, for each $1 \le j \le k_1$. All of these values are distinct and are different from both $-1$ and $2$, due to Lemma \ref{abcd_lemma}. They also cannot be equal to $0$, since $\mathcal{F}(0) = -2$ and all the $2\cos\left(\dfrac{2j \pi}{k} \right)$ are greater than $-2$, for $1 \le j \le k_1$. Since $2\cos\left(\dfrac{2\pi}{k} \right)$ is an eigenvalue of $B$ whose multiplicity is $2$, for each $1 \le j \le k_1$, we conclude that the $4k_1$-element set
\begin{align*}
    Z_1 &= \{ \alpha(y), \beta(y), \gamma(y), \delta(y) \colon y = 2\cos\left(\dfrac{2j \pi}{k}\right), j=\overline{1, k_1} \}
\end{align*}
is composed of eigenvalues of $A$, each of which has a multiplicity of at least $2$.

Each $\lambda \in \{ \alpha(2), \beta(2), \gamma(2), \delta(2) \}$ is clearly not equal to $-1$, $2$ or $0$, and satisfies $\mathcal{F}(\lambda) = 2$. Since $2$ is an eigenvalue of $B$, this means that  the $4$-element set
\begin{align*}
    Z_2 = \{ \alpha(2), \beta(2), \gamma(2), \delta(2) \}
\end{align*}
is composed of eigenvalues of $A$.

It is clear that $Z_1 \cap Z_2 = \varnothing$ and $\lvert Z_1 \rvert = 4k_1 = 2k - 2, \lvert Z_2 \rvert = 4$. Taking into consideration that $A$ is of order $5k$ and has the eigenvalue $-1$ with a multiplicity of $k$, we conclude that all of its eigenvalues from the set $Z_1$ must have a multiplicity of exactly $2$ and all of its eigenvalues from the set $Z_2$ must be simple. Also, the matrix $A$ cannot have any additional eigenvalues other than the ones we have mentioned. Having determined the spectrum of $A$, we finally obtain the following formula
\begin{align*}
    \mathcal{E}(W\!gc_k) = k\lvert -1 \rvert &+ \lvert \alpha(2) \rvert + \lvert \beta(2) \rvert + \lvert \gamma(2) \rvert + \lvert \delta(2) \rvert\\
    &+ 2 \sum_{j=1}^{k_1} \left| \alpha\left(2\cos\left(\dfrac{2j \pi}{k}\right)\right) \right| + 2\sum_{j=1}^{k_1} \left| \beta\left(2\cos\left(\dfrac{2j \pi}{k}\right)\right) \right|\\
    &+ 2\sum_{j=1}^{k_1} \left| \gamma\left(2\cos\left(\dfrac{2j \pi}{k}\right)\right) \right| + 2\sum_{j=1}^{k_1} \left| \delta\left(2\cos\left(\dfrac{2j \pi}{k}\right)\right) \right| .
\end{align*}
Bearing in mind that $\cos\left(\dfrac{2j \pi}{k}\right) = \cos\left(\dfrac{2(k - j) \pi}{k}\right)$ for $1 \le j \le k_1$, as well as $\cos\left(\dfrac{2 \cdot 0 \pi}{k}\right) = 2$, we promptly get
\begin{align}\label{wgc_sol16}
    \begin{split}
    \mathcal{E}(W\!gc_k) = k\lvert -1 \rvert &+ \sum_{j=0}^{k-1} \left| \alpha\left(2\cos\left(\dfrac{2j \pi}{k}\right)\right) \right| + \sum_{j=0}^{k-1} \left| \beta\left(2\cos\left(\dfrac{2j \pi}{k}\right)\right) \right|\\
    &+ \sum_{j=0}^{k-1} \left| \gamma\left(2\cos\left(\dfrac{2j \pi}{k}\right)\right) \right| + \sum_{j=0}^{k-1} \left| \delta\left(2\cos\left(\dfrac{2j \pi}{k}\right)\right) \right| .
    \end{split}
\end{align}

\noindent
\emph{Subcase }$k$ is even.\quad Let $k = 2k_1$. In this scenario, $B$ has two simple eigenvalues $2$ and $-2$, and $k_1 - 1$ eigenvalues $2\cos\left(\dfrac{2\pi}{k} \right), 2\cos\left(\dfrac{4\pi}{k} \right), \ldots, 2\cos\left(\dfrac{2(k_1-1) \pi}{k} \right)$ which all have a multiplicity of $2$.

Again, we know that the value $2\cos\left(\dfrac{2j \pi}{k} \right)$ is an element of the interval $(-2, 2)$ for each $1 \le j \le k_1-1$. Lemma \ref{abcd_lemma} implies the existence of exactly $4$ distinct real numbers $\alpha\left(2\cos\left(\dfrac{2j \pi}{k}\right)\right), \beta\left(2\cos\left(\dfrac{2j \pi}{k}\right)\right), \gamma\left(2\cos\left(\dfrac{2j \pi}{k}\right)\right), \linebreak \delta\left(2\cos\left(\dfrac{2j \pi}{k}\right)\right)$ that represent the solutions to the equation $\mathcal{F}(x) = 2\cos\left(\dfrac{2j \pi}{k}\right)$ in $x \in \mathbb{R}$, for all $1 \le j \le k_1-1$. It can be seen that all of these $4(k_1-1)$ values are distinct and they are all different from $-1$, $2$ and $0$, with a similar argumentation as the one used in the previous subcase. Due to the fact that $2\cos\left(\dfrac{2\pi}{k} \right)$ is an eigenvalue of $B$ with a multiplicity of $2$ for all $1 \le j \le k_1-1$, we obtain that the $(4k_1-4)$-element set
\begin{align*}
    Z_1 &= \{ \alpha(y), \beta(y), \gamma(y), \delta(y) \colon y = 2\cos\left(\dfrac{2j \pi}{k}\right), j=\overline{1, k_1-1} \}
\end{align*}
is composed of eigenvalues of $A$ whose multiplicities are at least $2$.

It is trivial to check that each $\lambda \in \{ \alpha(2), \beta(2), \gamma(2), \delta(2) \} \cup \{\alpha(-2), \gamma(-2), \linebreak \delta(-2) \}$ is not equal to $-1$, $2$ or $0$. These values also satisfy $\mathcal{F}(\lambda) \in \{2, -2\}$. Since $2$ and $-2$ are both eigenvalues of $B$, we get that the $7$-element set
$$Z_2 = \{ \alpha(2), \beta(2), \gamma(2), \delta(2), \alpha(-2), \gamma(-2), \delta(-2) \}$$
is composed of eigenvalues of $A$.

It is easy to see that $Z_1 \cap Z_2 = \varnothing$ and $\lvert Z_1 \rvert = 4k_1 - 4 = 2k-4, \lvert Z_2 \rvert = 7$. Since $A$ is of order $5k$ and we already know that it has an eigenvalue $-1$ of multiplicity $k$ and a simple eigenvalue $0$, we conclude that all of its eigenvalues from the set $Z_1$ must have a multiplicity of exactly $2$ and all of its eigenvalues from the set $Z_2$ must be simple. The matrix cannot have any additional eigenvalues other than the ones we have found, as well. Having found the entire spectrum of $A$, we reach the expression
\begin{align*}
    \mathcal{E}(W\!gc_k) = k\lvert -1 \rvert &+ \lvert \alpha(2) \rvert + \lvert \beta(2) \rvert + \lvert \gamma(2) \rvert + \lvert \delta(2) \rvert\\
    &+ \lvert \alpha(-2) \rvert + \lvert \gamma(-2) \rvert + \lvert \delta(-2) \rvert\\
    &+ 2 \sum_{j=1}^{k_1-1} \left| \alpha\left(2\cos\left(\dfrac{2j \pi}{k}\right)\right) \right| + 2\sum_{j=1}^{k_1-1} \left| \beta\left(2\cos\left(\dfrac{2j \pi}{k}\right)\right) \right|\\
    &+ 2\sum_{j=1}^{k_1-1} \left| \gamma\left(2\cos\left(\dfrac{2j \pi}{k}\right)\right) \right| + 2\sum_{j=1}^{k_1-1} \left| \delta\left(2\cos\left(\dfrac{2j \pi}{k}\right)\right) \right| .
\end{align*}
Taking into consideration that $\cos\left(\dfrac{2j \pi}{k}\right) = \cos\left(\dfrac{2(k - j) \pi}{k}\right)$ for $1 \le j \le k_1-1$, as well as $\cos\left(\dfrac{2 \cdot 0 \pi}{k}\right) = 2, \cos\left(\dfrac{2k_1 \pi}{k}\right) = -2$ and $\beta(-2) = 0$, we obtain the same Eq.\ (\ref{wgc_sol16}) as we did in the previous subcase. This allows us to finish the computation without having to divide the problem into multiple cases again.

Given the fact that $\alpha(y) < \beta(y) \le 0$ and $0 < \gamma(y) < \delta(y)$ for all $-2 \le y \le 2$, Eq.\ (\ref{wgc_sol16}) quickly transforms into
\begin{align*}
    \mathcal{E}(W\!gc_k) = k &- \sum_{j=0}^{k-1} \alpha\left(2\cos\left(\dfrac{2j \pi}{k}\right)\right) - \sum_{j=0}^{k-1} \beta\left(2\cos\left(\dfrac{2j \pi}{k}\right)\right)\\
    &+ \sum_{j=0}^{k-1} \gamma\left(2\cos\left(\dfrac{2j \pi}{k}\right)\right) + \sum_{j=0}^{k-1} \delta\left(2\cos\left(\dfrac{2j \pi}{k}\right)\right) .
\end{align*}
We are now in position to use Lemma \ref{abcd_sum} in order to get
\begin{alignat*}{2}
    && \alpha(y) + \beta(y) + \gamma(y) + \delta(y) &= 1\\
    \implies \quad && -\alpha(y)-\beta(y)+\gamma(y)+\delta(y) &= 1 - 2\alpha(y) - 2\beta(y) \, ,
\end{alignat*}
for each $-2 \le y \le 2$, which helps us finish the computation
\begin{align*}
    \mathcal{E}(W\!gc_k) &= k + \sum_{j=0}^{k-1} 1 - \sum_{j=0}^{k-1} 2 \alpha\left(2\cos\left(\dfrac{2j \pi}{k}\right)\right) - \sum_{j=0}^{k-1} 2 \beta\left(2\cos\left(\dfrac{2j \pi}{k}\right)\right)\\
    &= 2k - 2 \sum_{j=0}^{k-1} \alpha\left(2\cos\left(\dfrac{2j \pi}{k}\right)\right) - 2 \sum_{j=0}^{k-1} \beta\left(2\cos\left(\dfrac{2j \pi}{k}\right)\right). \qed
\end{align*}

The upcoming proof of Theorem \ref{main_theorem_3} will rely on a simple lemma from mathematical analysis. We present and prove this lemma before continuing with the main proof.
\begin{lemma}\label{integral_lemma}
    Let $\phi \colon [a, b] \to [\phi(a), \phi(b)]$ be a strictly monotonous and continuous function which is a bijection from $[a, b]$ to $[\phi(a), \phi(b)]$, and let $\phi^{-1} \colon [\phi(a), \phi(b)] \to [a, b]$ be its corresponding inverse function. We then have
    \begin{align}\label{integral_formula}
        \int_{a}^{b}\phi(x) \diff x + \int_{\phi(a)}^{\phi(b)} \phi^{-1}(x) \diff x &= b \ \phi(b) - a \ \phi (a) .
    \end{align}
\end{lemma}
\begin{proof}
    First of all, it is clear that $\phi^{-1}$ is also strictly monotonous and continuous. Thus, $\phi$ is Riemann-integrable on $[a, b]$ and $\phi^{-1}$ is Riemann-integrable on $[\phi(a), \phi(b)]$, which means that the integrals from Eq.\ (\ref{integral_formula}) are well defined.
    
    Let $P$ be a partition of $[a, b]$ composed of the points $c_0, c_1, c_2, \ldots, c_g$ which form a strictly monotonous sequence such that $c_0 = a$ and $c_g = b$. It is obvious that $\phi(c_0), \phi(c_1), \phi(c_2), \ldots, \phi(c_g)$ will be a strictly monotonous sequence of points such that $\phi(c_0) = \phi(a)$ and $\phi(c_g) = \phi(b)$, hence it represents a partition of $[\phi(a), \phi(b)]$, which we will denote via $Q$. Suppose that the distinguished points $\xi$ from $P$ are $c_0, c_1, c_2, \ldots, c_{g-1}$ respectively. Likewise, suppose that the distinguished points $\psi$ from $Q$ are $\phi(c_1), \phi(c_2), \phi(c_3), \ldots, \phi(c_{g})$ respectively. This leads us to the following two Riemann sums
    \begin{align*}
        \sigma(\phi, P, \xi) &= \sum_{j = 0}^{g-1} \phi(c_j) \cdot (c_{j+1} - c_j) \, ,\\
        \sigma(\phi^{-1}, Q, \psi) &= \sum_{j = 0}^{g-1} c_{j+1} \cdot (\phi(c_{j+1}) - \phi(c_j)) \, .
    \end{align*}
    From the definition of the Riemann integral, we know that
    \begin{align*}
        \lim_{\lambda(P) \to 0} \sigma(\phi, P, \xi) &= \int_{a}^{b}\phi(x) \diff x \, .
    \end{align*}
    Due to the fact that $\phi$ is continuous on the compact set $[a, b]$, it follows that $\phi$ must be uniformly continuous on $[a, b]$ as well. This means that $\lambda(P) \to 0$ implies $\lambda(Q) \to 0$, which leads us to
    \begin{alignat*}{2}
        && \lim_{\lambda(Q) \to 0} \sigma(\phi^{-1}, Q, \psi) &= \int_{\phi(a)}^{\phi(b)} \phi^{-1}(x) \diff x\\
        \implies \quad && \lim_{\lambda(P) \to 0} \sigma(\phi^{-1}, Q, \psi) &= \int_{\phi(a)}^{\phi(b)} \phi^{-1}(x) \diff x \, .
    \end{alignat*}
    Thus, we get
    \begin{align*}
        \int_{a}^{b}\phi(x) \diff x + \int_{\phi(a)}^{\phi(b)} \phi^{-1}(x) \diff x &= \lim_{\lambda(P) \to 0} \left( \sigma(\phi, P, \xi) + \sigma(\phi^{-1}, Q, \psi) \right) .
    \end{align*}
    However, we know that
    \begin{align*}
        \sigma(\phi, P, \xi) + \sigma(\phi^{-1}, Q, \psi) &= \sum_{j = 0}^{g-1} \phi(c_j) \cdot (c_{j+1} - c_j) + \sum_{j = 0}^{g-1} c_{j+1} \cdot (\phi(c_{j+1}) - \phi(c_j))\\
        &= \sum_{j=0}^{g-1} \left( \phi(c_j) c_{j+1} - \phi(c_j) c_j + c_{j+1} \phi(c_{j+1}) - c_{j+1} \phi(c_j) \right)\\
        &= \sum_{j=0}^{g-1} \left( c_{j+1} \phi(c_{j+1}) - c_j \phi(c_j) \right)\\
        &= c_g \phi(c_g) - c_0 \phi(c_0)\\
        &= b \ \phi(b) - a \ \phi(a) \, ,
    \end{align*}
    which directly implies Eq.\ (\ref{integral_formula}).
\end{proof}

\bigskip\noindent
{\em Proof of Theorem \ref{main_theorem_3}}.\quad The functions $\alpha$ and $\beta$ are both continuous on their entire domain $[-2, 2]$, due to Lemma \ref{abcd_lemma}. It is clear that the function $x \to 2 \cos x$ is defined and continuous on $\mathbb{R}$, with values in $[-2, 2]$. Thus, the respective composite functions $x \to \alpha(2 \cos x)$ and $x \to \beta(2 \cos x)$ must also be defined and continuous on $\mathbb{R}$. This implies that these functions are Riemann-integrable on each bounded and closed interval.

Taking into consideration that $\mu(W\!gp_k) = 2k$, Eq.\ (\ref{main_formula_1}) gives us
\begin{align*}
    \dfrac{\mathcal{E}(W\!gp_k)}{\mu(W\!gp_k)} &= 1 - \dfrac{1}{k} \sum_{j=1}^{k} \alpha \left( 2 \cos\left( \frac{j \pi}{k+1}\right) \right) - \dfrac{1}{k} \sum_{j=1}^{k} \beta\left( 2 \cos\left( \frac{j \pi}{k+1}\right) \right) .
\end{align*}
Here, it is important to notice that
\begin{align*}
    \dfrac{1}{k} \sum_{j=1}^{k} \alpha \left( 2 \cos\left( \frac{j \pi}{k+1}\right) \right) &= \dfrac{1}{k} \left( -\alpha\left( 2 \cos\left( \frac{0 \pi}{k+1}\right)\right) + \sum_{j=0}^{k} \alpha \left( 2 \cos\left( \frac{j \pi}{k+1}\right) \right) \right)\\
    &= - \dfrac{\alpha(2)}{k} + \dfrac{1}{k} \sum_{j=0}^{k} \alpha \left( 2 \cos\left( \frac{j \pi}{k+1}\right) \right)\\
    &= - \dfrac{\alpha(2)}{k} + \dfrac{k+1}{k \pi} \sum_{j=0}^{k} \alpha \left( 2 \cos\left( \frac{j \pi}{k+1}\right) \right) \cdot \dfrac{\pi}{k+1} \, ,
\end{align*}
where $\displaystyle\sum_{j=0}^{k} \alpha \left( 2 \cos\left( \dfrac{j \pi}{k+1}\right) \right) \cdot \dfrac{\pi}{k+1}$ actually represents a Riemann sum of the function $x \to \alpha(2 \cos x)$ over the closed interval $[0, \pi]$. The mesh of this Riemann sum obviously equals $\dfrac{\pi}{k+1}$, which directly implies
\begin{align*}
    \lim_{k \to \infty} \sum_{j=0}^{k} \alpha \left( 2 \cos\left( \frac{j \pi}{k+1}\right) \right) \cdot \dfrac{\pi}{k+1} &= \int_{0}^{\pi}\alpha(2\cos x) \diff x \, .
\end{align*}
Thus, it is straightforward to see that
\begin{align*}
    \lim_{k \to \infty} \dfrac{1}{k} \sum_{j=1}^{k} \alpha \left( 2 \cos\left( \frac{j \pi}{k+1}\right) \right) &= \dfrac{1}{\pi} \int_{0}^{\pi}\alpha(2\cos x) \diff x \, .
\end{align*}
In a completely analogous manner, it can be shown that
\begin{align*}
    \lim_{k \to \infty} \dfrac{1}{k} \sum_{j=1}^{k} \beta \left( 2 \cos\left( \frac{j \pi}{k+1}\right) \right) &= \dfrac{1}{\pi} \int_{0}^{\pi}\beta(2\cos x) \diff x \, ,
\end{align*}
which leads us to
\begin{align}\label{wgp_almost_done}
    \lim_{k \to \infty} \dfrac{\mathcal{E}(W\!gp_k)}{\mu(W\!gp_k)} &= 1 - \dfrac{1}{\pi} \int_{0}^{\pi}\alpha(2\cos x) \diff x - \dfrac{1}{\pi} \int_{0}^{\pi}\beta(2\cos x) \diff x \, .
\end{align}
On the other hand, we know that $\mu(W\!gc_k) = 2k$. Eq.\ (\ref{main_formula_2}) implies
\begin{align*}
    \dfrac{\mathcal{E}(W\!gc_k)}{\mu(W\!gc_k)} = 1 - \dfrac{1}{k} \sum_{j=0}^{k-1} \alpha \left( 2 \cos\left( \frac{2j \pi}{k}\right) \right) - \dfrac{1}{k} \sum_{j=0}^{k-1} \beta\left( 2 \cos\left( \frac{2j \pi}{k}\right) \right) .
\end{align*}
Again, we notice that
\begin{align*}
    \dfrac{1}{k} \sum_{j=0}^{k-1} \alpha \left( 2 \cos\left( \frac{2j \pi}{k}\right) \right) &= \dfrac{1}{2 \pi} \sum_{j=0}^{k-1} \alpha \left( 2 \cos\left( \frac{2j \pi}{k}\right) \right) \cdot \dfrac{2 \pi}{k} \, ,
\end{align*}
where $\displaystyle\sum_{j=0}^{k-1} \alpha \left( 2 \cos\left( \dfrac{2j \pi}{k}\right) \right) \cdot \dfrac{2 \pi}{k}$ is a Riemann sum of the function $x \to \alpha(2 \cos x)$ over the closed interval $[0, 2\pi]$. The mesh of this Riemann sum equals $\dfrac{2 \pi}{k}$, hence we get
\begin{align*}
    \lim_{k \to \infty} \sum_{j=0}^{k-1} \alpha \left( 2 \cos\left( \frac{2j \pi}{k}\right) \right) \cdot \dfrac{2 \pi}{k} &= \int_{0}^{2 \pi} \alpha(2 \cos x) \diff x \, ,
\end{align*}
along with
\begin{align*}
    \lim_{k \to \infty} \dfrac{1}{k} \sum_{j=0}^{k-1} \alpha \left( 2 \cos\left( \frac{2j \pi}{k}\right) \right) &= \dfrac{1}{2 \pi} \int_{0}^{2 \pi} \alpha(2 \cos x) \diff x \, .
\end{align*}
The same can analogously be proven for the function $x \to \beta(2 \cos x)$, i.e.\
\begin{align*}
    \lim_{k \to \infty} \dfrac{1}{k} \sum_{j=0}^{k-1} \beta \left( 2 \cos\left( \frac{2j \pi}{k}\right) \right) &= \dfrac{1}{2 \pi} \int_{0}^{2 \pi} \beta(2 \cos x) \diff x \, .
\end{align*}
This leads us to
\begin{align*}
    \lim_{k \to \infty} \dfrac{\mathcal{E}(W\!gc_k)}{\mu(W\!gc_k)} &= 1 - \dfrac{1}{2 \pi} \int_{0}^{2 \pi} \alpha(2 \cos x) \diff x - \dfrac{1}{2 \pi} \int_{0}^{2 \pi} \beta(2 \cos x) \diff x \, .
\end{align*}
Here, it is useful to notice that via substitution $x = 2\pi - z$, we get
\begin{align*}
    \int_{\pi}^{2 \pi} \alpha(2 \cos x) \diff x &= \int_{\pi}^{0} \alpha(2 \cos (2\pi - z)) (-1) \diff z\\
    &= \int_{0}^{\pi} \alpha(2 \cos z) \diff z \, ,
\end{align*}
which means that
\begin{align*}
    \int_{0}^{2 \pi} \alpha(2 \cos x) \diff x &= \int_{0}^{\pi} \alpha(2 \cos x) \diff x + \int_{\pi}^{2 \pi} \alpha(2 \cos x) \diff x\\
    &= 2  \int_{0}^{\pi} \alpha(2 \cos x) \diff x \, .
\end{align*}
Similarly, it can be shown that $\displaystyle\int_{0}^{2 \pi} \beta(2 \cos x) \diff x = 2 \displaystyle\int_{0}^{\pi} \beta(2 \cos x) \diff x$, which helps us obtain
\begin{align}\label{wgc_almost_done}
    \lim_{k \to \infty} \dfrac{\mathcal{E}(W\!gc_k)}{\mu(W\!gc_k)} &= 1 - \dfrac{1}{\pi} \int_{0}^{\pi} \alpha(2 \cos x) \diff x - \dfrac{1}{\pi} \int_{0}^{\pi} \beta(2 \cos x) \diff x \, .
\end{align}
By comparing Eq.\ (\ref{wgc_almost_done}) to Eq.\ (\ref{wgp_almost_done}), we see that $\displaystyle\lim_{k \to \infty} \dfrac{\mathcal{E}(W\!gp_k)}{\mu(W\!gp_k)} = \displaystyle\lim_{k \to \infty} \dfrac{\mathcal{E}(W\!gc_k)}{\mu(W\!gc_k)}$. We will denote the value of these two limits by $L$. In order to complete the proof, it is sufficient to show that this value of $L$ satisfies Eq.\ (\ref{main_formula_3}).

The function $x \to \alpha(2 \cos x)$ is continuous and strictly increasing on $[0, \pi]$. It also represents a bijection from $[0, \pi]$ to $[\alpha(2), \alpha(-2)]$, with the corresponding inverse function having the form $x \to \arccos\left(\dfrac{\mathcal{F}(x)}{2}\right)$. This allows us to apply Lemma \ref{integral_lemma} in order to obtain
\begin{alignat*}{2}
    && \int_{0}^{\pi}\alpha(2 \cos x) \diff x + \int_{\alpha(2)}^{\alpha(-2)} \arccos\left(\dfrac{\mathcal{F}(x)}{2}\right) \diff x &= \alpha(-2) \cdot \pi - \alpha(2) \cdot 0\\
    \implies \quad && \int_{0}^{\pi}\alpha(2 \cos x) \diff x + \int_{\alpha(2)}^{\alpha(-2)} \arccos\left(\dfrac{\mathcal{F}(x)}{2}\right) \diff x &= \pi \ \alpha(-2)  \, .
\end{alignat*}
The function $x \to \beta(2 \cos x)$ is also continuous and strictly increasing on $[0, \pi]$. It is easy to see that this function represents a bijection from $[0, \pi]$ to $[\beta(2), \beta(-2)]$, with the corresponding inverse function being equal to $x \to \arccos\left(\dfrac{\mathcal{F}(x)}{2}\right)$. By implementing Lemma \ref{integral_lemma}, we conclude that
\begin{alignat*}{2}
    && \int_{0}^{\pi}\beta(2 \cos x) \diff x + \int_{\beta(2)}^{\beta(-2)} \arccos\left(\dfrac{\mathcal{F}(x)}{2}\right) \diff x &= \beta(-2) \cdot \pi - \beta(2) \cdot 0\\
    \implies \quad && \int_{0}^{\pi}\beta(2 \cos x) \diff x + \int_{\beta(2)}^{0} \arccos\left(\dfrac{\mathcal{F}(x)}{2}\right) \diff x &= 0 \, ,
\end{alignat*}
due to the fact that $\beta(-2) = 0$. Thus, we have
\begin{align*}
    \dfrac{1}{\pi} \int_{0}^{\pi}\alpha(2 \cos x) \diff x &= \alpha(-2) - \dfrac{1}{\pi} \int_{\alpha(2)}^{\alpha(-2)} \arccos\left(\dfrac{\mathcal{F}(x)}{2}\right) \diff x \, ,\\
    \dfrac{1}{\pi} \int_{0}^{\pi}\beta(2 \cos x) \diff x &= -\dfrac{1}{\pi} \int_{\beta(2)}^{0} \arccos\left(\dfrac{\mathcal{F}(x)}{2}\right) \diff x \, .
\end{align*}
Hence, we finally obtain
\begin{align*}
    L &= 1 - \dfrac{1}{\pi} \int_{0}^{\pi} \alpha(2 \cos x) \diff x - \dfrac{1}{\pi} \int_{0}^{\pi} \beta(2 \cos x) \diff x\\
    &= 1 - \alpha(-2) + \dfrac{1}{\pi} \int_{\alpha(2)}^{\alpha(-2)} \arccos\left(\dfrac{\mathcal{F}(x)}{2}\right) \diff x + \dfrac{1}{\pi} \int_{\beta(2)}^{0} \arccos\left(\dfrac{\mathcal{F}(x)}{2}\right) \diff x \, .
\end{align*}
which completes the proof.\qed

\bigskip
It is straightforward to numerically check that $L > 2\sqrt{3}$. Moreover, it can be computed that $L \approx 3.483650329$ with a precision of $10$ digits, while $2 \sqrt{3} \approx 3.464101615$. This means that only finitely many elements $G$ of the infinite graph sequences $(W\!gp_k)_{k \in \mathbb{N}}$ and $(W\!gc_k)_{k \in \mathbb{N} \setminus \{ 1 \} }$ can satisfy the inequality $\mathcal{E}(G) \le 2 \mu(G) \sqrt{\Delta}$. Hence, we have proved that the two graph sequences both yield an infinite family of counterexamples to the inequality of interest.


\vspace{1cc}


{\small
\noindent
{\bf \DJ{}or\dj{}e Stevanovi\'c,} \\
Ktitor 35, \\
Ni\v s, Serbia, \\
e-mail: st.djole@yahoo.com

\medskip\noindent
{\bf Ivan Damnjanovi\'c,} \\
University of Ni\v s, Faculty of Electronics, \\
Ni\v s, Serbia, \\
e-mail: ivan.damnjanovic@elfak.ni.ac.rs

\medskip\noindent
{\bf Dragan Stevanovi\'c,} \\
Mathematical Institute of the Serbian Academy of Sciences and Arts, \\
Belgrade, Serbia, \\
e-mail: dragan\_stevanovic@mi.sanu.ac.rs
}

\end{document}